\documentclass[article,reqno]{amsart}
\usepackage[utf8]{inputenc}
\usepackage{a4wide, amsfonts}
\usepackage[english]{babel}
\usepackage{xcolor}
\usepackage[T1]{fontenc}
\usepackage{amsmath, amssymb,epic,graphicx,mathrsfs,enumerate}    
\usepackage{cite}
\usepackage[numbers,sort&compress]{natbib}
\usepackage[all]{xy}
\usepackage{color}
\usepackage{comment}
\usepackage{enumitem}
\usepackage{esint}
\usepackage{hyperref}
\usepackage{amsthm}
\usepackage{latexsym}
\usepackage{epsfig}
\usepackage{soul}
\usepackage{geometry}
\usepackage{doi}
\hypersetup{colorlinks=true, linkcolor=black, filecolor=black, urlcolor=blue,
citecolor=blue}
\newcommand{\blu}[1]{{\color{blue}#1}}
\newcommand{\red}[1]{{\color{red}#1}}

\newtheorem{theorem}{Theorem}[section]
\newtheorem{corollary}[theorem]{Corollary}
\newtheorem{lemma}[theorem]{Lemma}
\newtheorem{prop}[theorem]{Proposition}

\theoremstyle{remark}
\newtheorem{obs}{Remark}

\newcommand{\R}{\mathbb{R}}
\newcommand{\K}{\mathcal{K}}
\newcommand{\E}{\mathcal{E}}
\newcommand{\e}{\varepsilon}
\newcommand{\al}{\alpha}
\newcommand{\g}{\gamma}
\newcommand{\la}{\lambda}
\newcommand{\s}{\sigma}
\newcommand{\p}{\varphi}
\newcommand{\de}{\delta}

\title[Concentration for Choquard MFG]{Mass concentration for Ergodic Choquard Mean-Field Games}

\author[C.Bernardini]{Chiara Bernardini}
\address{Dipartimento di Matematica "Tullio Levi Civita",
Università di Padova, Via Trieste 63, 35121 Padova, Italy}
\email{chiara.bernardini@math.unipd.it}
 
\begin{document}

\begin{abstract}
We study concentration phenomena in the vanishing viscosity limit for second-order stationary Mean-Field Games systems defined in the whole space $\R^N$ with Riesz-type aggregating coupling and external confining potential. In this setting, every player of the game is attracted toward congested areas and the external potential discourages agents from being far away from the origin. Assuming some restrictions on the strength of the attractive nonlocal term depending on the growth of the Hamiltonian, we study the asymptotic behavior of solutions in the vanishing viscosity limit. First, we obtain existence of classical solutions to potential-free MFG systems with Riesz-type coupling. Secondly, we prove concentration of mass around minima of the potential.\\ 


\end{abstract}

\maketitle

\date{}


\section{Introduction}

In the present work, we keep on the study, initiated in \cite{BerCes}, of stationary Mean-Field Games systems defined in the whole space $\R^N$ with aggregating nonlocal coupling of Riesz-type. In particular, given $M>0$, we consider systems of the form
\begin{equation}\label{1 senza V}
\begin{cases}
-\Delta u+\frac{1}{\g}|\nabla u|^\g+\la=-K_\al\ast m(x)\\
-\Delta m-\mathrm{div}(m\nabla u|\nabla u|^{\g-2})=0\\
\int_{\R^N}m=M,\quad m\ge0
\end{cases}\text{in}\,\,\R^N
\end{equation}

\noindent where $\g>1$ is fixed and $K_\al:\R^N\to\R$ is the Riesz potential of order $\al\in(0,N)$, which is defined for every $x\in\R^N\setminus\{0\}$ by 
$$K_\al(x)=\frac{1}{|x|^{N-\al}}.$$

\noindent Such nonlinear elliptic systems of PDEs arise in the context of Mean-Field Games (MFG in short) to describe Nash equilibria of differential games with infinitely many interacting players. We refer the reader to \cite{LL1,CaPo} and the references therein, for more details on the general theory of Mean-Field Games and its developments. In our setting, the coupling between the individual and the overall population is defined in terms of a Riesz-type interaction kernel $-K_\al\ast m$, since it is attractive and nonlocal, each player of the game is attracted toward regions where the population is highly distributed. Moreover, every player is subject to a Brownian noise $\sqrt{2}B_t$ which induces a dissipation effect. Existence results for classical solutions to the MFG system \eqref{1 senza V} will depend on a balancing between dissipation and aggregation. Indeed, if the long-range attractive force is too strong, the mass $m$ tends to concentrate and hence to develop singularities, while if the diffusion dominates, we might have loss of mass at infinity, in both cases we expect nonexistence of classical solutions. 

Stationary \textit{focusing} MFG systems with different assumptions on the coupling have been studied e.g. in \cite{CC,CeCi,Cir,Cir1,GNP}. The problem of existence of solutions to focusing MFG systems requires different approaches than the ones that have been developed in the literature to study \textit{defocusing} MFG systems, namely models where individuals avoid areas with high density of population. Indeed, many existence and regularity arguments and also uniqueness of equilibria requires to have an increasing coupling. In addition, our model is defined in the whole space $\R^N$ (refer e.g. to \cite{BaPi,BerCes,CC,GoPi,kouhkouh} for other works in the non-compact setting), this non-compact situation  brings new difficulties, due to the fact that the Brownian noise has to be balanced by the optimal velocity to avoid loss of mass. Finally, if the coupling is given by $K_\al\ast m$, dissipating forces dominate and we expect nonexistence of solutions in the noncompact setting $\R^N$; while if we consider a bounded domain with periodic or Neumann boundary conditions a similar approach to the case of power-like coupling should work (see e.g. \cite{Cir,GoPaVo}).
 
In this paper, focusing on the \textit{mass-subcritical regime} $N-\g'<\al<N$, where $\g'=\frac{\g}{\g-1}$ is the conjugate exponent of $\g$, we provide existence of classical solutions to the MFG system \eqref{1 senza V}. Notice that by \textit{classical solution} we mean a triple $(u,m,\la)\in C^2(\R^N)\times W^{1,p}(\R^N)\times \R$ for every $p\in(1,+\infty)$, solving the system. More precisely, we obtain the following existence result.

\begin{theorem}\label{teo no pot}
Let $N-\g'<\al<N$. Then, for every $M>0$ there exists $(\bar{u},\bar{m},\bar{\la})$ classical solution to the MFG system \eqref{1 senza V}. Moreover, there exist $C_1, C_2, C_3$ and $C_4$ positive constants such that 
$$\bar{u}(x)\ge C_1|x|-C_1^{-1},\quad\quad\quad |\nabla\bar{u}|\le C_2$$
and
$$0<\bar{m}(x)\le C_3e^{-C_4|x|}.$$
\end{theorem}

\begin{obs} 
Solutions to the MFG system \eqref{1 senza V} are invariant by translation, namely if $(\bar{u}(x),\bar{m}(x),\bar{\la})$ is a classical solution to \eqref{1 senza V} then for every $x_0\in\R^N$ and $c\in\R$, also $\left(\bar{u}(x+x_0)+c,\bar{m}(x+x_0),\bar{\la} \right)$ is a classical solution to \eqref{1 senza V}. Therefore, the constants $C_1$ and $C_4$ appearing in the previous theorem, depend on the choice of the solution.
\end{obs}

Theorem \ref{teo no pot} partially completes the study of existence of solutions to the potential-free MFG system \eqref{1 senza V} started in \cite{BerCes}. In particular, in \cite[Theorem 1.1]{BerCes}, using a Pohozaev-type identity, one proves that if $0<\al<N-2\g'$, ``regular'' solutions to the MFG system \eqref{1 senza V} (namely satisfying some quite natural integrability conditions and boundary conditions at infinity) do not exist. It remains still open the problem of existence of solutions to \eqref{1 senza V} when $\al\in[N-2\g',N-\g']$.

On the other hand, the main result in \cite{BerCes} deals with the study of MFG systems with Riesz-type coupling and external confining potential $V$. More in detail, exploiting a Schauder fixed point argument, one proves that if $\al\in(N-\g',N)$ the MFG system admits a classical solution for every total mass $M>0$, while if $\al\in(N-2\g',N-\g']$ a solution does exist at least for sufficiently small masses $M$, below some threshold value $M_0$ (see \cite[Theorem 1.2]{BerCes} for more details).
Notice that, using a fixed point approach, the presence of the coercive potential $V$ adds compactness to the problem and proves to be essential to conclude, so the existence result in \cite{BerCes} does not cover the case when the potential $V$ is identically $0$.
In order to deal with the potential-free system we take advantage a variational argument. This approach allows us to obtain some uniform (namely not depending on the viscosity parameter) estimates on the solutions, which will be crucial in the vanishing viscosity setting and which can not be obtained by means of a fixed point technique.

Finally, a similar MFG system but with local decreasing coupling defined in terms of a power-type function, has been studied in \cite{CC}. We point out that, in our setting the nonlocal attractive coupling models a long-range attractive force between players, moreover, in order to deal with the Riesz-term we need different techniques (see Section 2) compared to the ones used in \cite{CC}.

\begin{obs}
In the case $\g=\g'=2$, using the so called Hopf-Cole transformation $v(x):= e^{-u(x)/2}$  and setting $m(x)=v^2(x)$, we can reduce the MFG system \eqref{1 senza V} to the following normalized Choquard equation  
\begin{equation}\label{cho} 
\begin{cases}
-2\Delta v-\la v=(K_\al\ast v^2)v\\
\int_{\R^N}v^2(x) dx=M, \qquad v>0
\end{cases}\quad\text{in}\,\,\R^N.
\end{equation} 
Equation \eqref{cho} was first studied by E. Lieb \cite{Lieb0}, who proved existence and uniqueness (up to translations) of solutions when $N=3$ and $\al=2$ by using symmetric decreasing rearrangement inequalities. Then, P.-L. Lions \cite{Lions} proved that there exists a minimum of the energy associated to \eqref{cho} when we restrict the infimum to functions with spherical symmetry, refer to \cite[§3]{Li} and also \cite{MS,MVs} for further results. Here, from Theorem \ref{teo no pot} we obtain a more general result for the range of values $\al$ such that the normalized Choquard equation \eqref{cho} has a solution, but we left open the problem of symmetry of solutions.
\end{obs}

In order to study system \eqref{1 senza V}, we consider an ergodic MFG system defined in the whole space $\R^N$ with an external confining potential $V$ and Brownian noise which depends on $\e>0$. More in detail, we take into account systems of the form
\begin{equation}\label{1}
\begin{cases}
-\e\Delta u+\frac{1}{\g}|\nabla u|^\g+\la=V(x)-K_\al\ast m(x)\\
-\e\Delta m-\mathrm{div}(m\nabla u|\nabla u|^{\g-2})=0\\
\int_{\R^N}m=M,\quad m\ge0
\end{cases}\text{in}\,\,\R^N
\end{equation}
where we assume that the potential $V$ is locally H\"older continuous and there exist two positive constants $b$ and $C_V$ such that 
\begin{equation}\label{V}
C_V^{-1}(\max\{|x|-C_V,0\})^b\le V(x)\le C_V(1+|x|)^b,\quad \forall x\in\R^N.
\end{equation}
Studying the asymptotic behavior of rescaled solutions to the MFG system \eqref{1} in the vanishing viscosity limit, we are able to prove existence of classical solutions to the MFG system \eqref{1 senza V} (without the potential term $V$). As a matter of fact, letting $\e\to0$, the dynamic of each player is no subject anymore to the dissipation effect induced by the Brownian motion, so we expect aggregation of players. In particular, in the vanishing viscosity limit the mass $m$ tends to concentrate, the introduction of the coercive potential $V$, which represents spatial preferences of agents, rules out this possibility and leads to concentration of mass around minima of the potential $V$.

Let us summarize the main tools to prove our results. As J.-M. Lasry and P.-L. Lions first pointed out in \cite{LL3}, taking into account the variational nature of the MFG system, solutions to \eqref{1} are related to critical points of the following energy functional \begin{equation}\label{energy}
\mathcal{E}(m,w):=\begin{cases} 
\int\limits_{\R^N}mL\left(-\frac{w}{m}\right)+V(x)\,m \,dx-\frac{1}{2}\int\limits_ {\R^N}\int\limits_{\R^N}\frac{m(x)\,m(y)}{|x-y|^{N-\al}} dx\,dy\qquad \text{if}\,\,(m,w)\in\mathcal{K}_{\e,M},\\
+\infty \hspace{7.9cm} \text{otherwise}\end{cases}
\end{equation}
where
$$L\left(-\frac{w}{m}\right):=\begin{cases} 
\frac{1}{\g'}\left|\frac{w}{m} \right|^{\g'}\qquad\text{if } m>0\\
0\qquad\qquad\quad\,\text{if } m=0,\,w=0\\
+\infty\qquad\quad\,\,\,\,\,\text{otherwise}
\end{cases}$$
and the constraint set is defined as
\begin{equation}\label{K}
\begin{aligned}
\mathcal{K}_{\e,M}:=\Big\{(m,w)\in & (L^1(\R^N)\cap L^q(\R^N))\times L^1(\R^N)\quad\text{s.t.}\quad \int_{\R^N}m\,dx=M,\quad m\ge0\,\,\, \text{a.e.}\\
&\e\int_{\R^N}m(-\Delta\p)\,dx=\int_{\R^N}w\cdot\nabla\p\,dx\quad\forall\p\in C^\infty_0(\R^N)\Big\}
\end{aligned}
\end{equation}
with
\begin{equation}\label{q}
q:=\begin{cases}\frac{N}{N-\g'+1}\quad\text{if}\,\,\g'< N\\
\g'\qquad\quad\,\,\,\text{if}\,\,\g'\ge N\end{cases}.
\end{equation}
Using some regularity results for the Kolmogorov equation (refer to Proposition \ref{kolm} below), the Hardy-Littlewood-Sobolev inequality and the fact that $V$ is non-negative, we prove that the energy $\E$ is bounded from below. 
By classical direct methods and compactness arguments, we obtain minimizers $(m_\e,w_\e)$ of $\E$. Finally, passing to another functional with linearized Riesz-term and using convex duality arguments (see for instance \cite{BC,CC,CeCi,CG,CGPT}), we are able to construct the associated solutions $(u_\e,m_\e,\la_\e)$ of the MFG system \eqref{1}. Then, in order to investigate the behavior of the system in the vanishing viscosity limit, we define a suitable rescaling of $u_\e$, $m_\e$ and $\la_\e$. We also translate the reference system by $y_\e$, where $y_\e$ is a point of minimum for the value function $u_\e$, in this way around $y_\e$ the mass remains positive and we can rule out vanishing of the total mass in the limit. We obtain a triple $(\bar{m}_\e,\bar{u}_\e,\tilde{\la}_\e)$ which solves the following MFG system
$$\begin{cases}
-\Delta\bar{u}_\e+\frac{1}{\g}|\nabla\bar{u}_\e|^\g+\tilde{\la}_\e=\e^\frac{(N-\al)\g'}{\g'-N+\al}V\left(\e^\frac{\g'}{\g'-N+\al}(y+y_\e)\right)-K_\al\ast\bar{m}_\e(y)\\
-\Delta\bar{m}_\e-\mathrm{div}(\bar{m}_\e\nabla\bar{u}_\e|\nabla\bar{u}_\e|^{\g-2})=0\\
\int_{\R^N}\bar{m}_\e=M
\end{cases}$$
(see Subsection \ref{pb riscalato} for more details). Exploiting a concentration-compactness argument (refer to the seminal work of P.-L. Lions \cite{Li}) as done in \cite{CC}, we are able to prove that there is no loss of mass when passing to the limit as $\e\to 0$. We show that in the vanishing viscosity limit, the rescaled solutions converge (up to sub-sequences) to $(\bar{u},\bar{m},\bar{\la})$ \textit{classical} solution to the MFG system \eqref{1 senza V}. Moreover, solutions to \eqref{1 senza V} are related to minimum points of the following energy
$$\E_0(m,w):=\int_{\R^N}\frac{m}{\g'}\bigg|\frac{w}{m}\bigg|^{\g'}dx-\frac{1}{2}\int_{\R^N}\int_{\R^N}\frac{m(x)m(y)}{|x-y|^{N-\al}}dx\,dy$$
over the constraint set
$$\mathcal{B}:=\Big\{(m,w)\in\K_{1,M}\,\,\Big|\,\,m(1+|x|^b)\in L^1(\R^N)\Big\}.$$
The following theorem states existence of solutions to \eqref{1} and concentration of mass.

\begin{theorem}\label{TEO concentrazione massa}
Let $N-\g'<\al<N$. Assume that the potential $V$ is locally H\"older continuous and satisfies \eqref{V}. Then, for every $\e,M>0$ there exists $(u_\e,m_\e,\la_\e)$ classical solution to \eqref{1}, such that $(m_\e,-m_\e\nabla u_\e|\nabla u_\e|^{\g-2})$ is a minimum of the energy $\E$. \\
\indent Moreover, there exists a sequence $\e\to0$ and a sequence of points $x_\e$ around which there is concentration of mass, namely for every $\eta>0$ there exist $R,\e_0>0$ such that 
$$\int_{B(x_\e,\e^\frac{\g'}{\g'-N+\al}R)}m_\e(x)\,dx\ge M-\eta$$
for all $\e<\e_0$ and 
$$x_\e\to\bar{x}, \quad\text{as}\,\,\,\e\to0$$
where $\bar{x}$ is a minimum point of the potential $V$ and $V(\bar{x})=0$. 
\end{theorem} 

\begin{obs}\label{condizioni Ham}
We assumed that the Hamiltonian $H$ has the form $H(p)=\frac{1}{\g} |p|^\g$ for $\g>1$ fixed, but actually the previous results hold also for more general assumptions on the Hamiltonian $H$, namely assuming that the Hamiltonian $H:\R^N\to\R$ is strictly convex, $H\in C^2(\R^N\setminus\{0\})$ and there exist $C_H,\,K>0$ and $\g>1$, such that $\forall p\in \R^N$ the following conditions hold
$$C_H|p|^\g-K\le H(p)\le C_H|p|^\g$$
$$\nabla H(p)\cdot p-H(p)\ge K^{-1}|p|^\g-K$$
$$|\nabla H(p)|\le K |p|^{\g-1}.$$
\end{obs}

The outline of the paper is the following. In Section 2 we provide some preliminary results. In particular, we recall some a priori estimates and elliptic regularity results for solutions to the Kolmogorov equation, and also integrability, H\"older continuity and compactness results for the Riesz potential term. In Section 3, using a variational approach, we prove existence of minimizers of the energy and from them we obtain the associated solution to the MFG system. In such a way we are in the position to analyze the asymptotic behavior of solutions in the vanishing viscosity limit, in particular, in Section 4 we prove existence of ground states to MFG system defined in the whole space $\R^N$ with Riesz-type coupling and without the confining potential $V$. Finally, in Section 5, we show concentration of the mass toward minima of $V$.

In what follows, $C, C_1, C_2, K_1,\dots$ denote generic positive constants which may change from line to line and also within the same line.\\


\section{Preliminaries}\label{sez2}

We recall here some a priori estimates and elliptic regularity results for solutions to the Kolmogorov equation; we mention also some regularity results for solutions to Hamilton-Jacobi-Bellman equations defined in the whole space $\R^N$. Finally, we state some properties of the Riesz potential. For further details we refer to \cite{BerCes,CC}. Let us fix $\e,M>0$. We will always assume that either $\g'\ge N$ or $\g'<N$ and $N-\g'\le\al<N$.

\subsection{Regularity for the Kolmogorov equation and the Hamilton-Jacobi-Bellman equation.}

We will use the following result (proved in \cite[Proposition 2.4]{CC}).

\begin{prop}\label{prop2.4CC}
Let $m\in L^p(\R^N)$ for $p>1$ and assume that for some $K>0$ 
$$\left|\int_{\R^N}m\Delta\varphi\,dx\right|\le K\|\nabla\varphi\|_{L^{p'} (\R^N)}, \quad \forall\varphi\in C_0^\infty(\R^N).$$
Then, $m\in W^{1,p}(\R^N)$ and there exists a constant $C>0$ depending only on $p$ such that 
$$\|\nabla m\|_{L^p(\R^N)}\le C\,K.$$
\end{prop}

\begin{prop}\label{kolm}
Assume that $(m,w)\in \K_{\e,M}$ and $E:=\frac{1}{\e^{\g'}}\int_{\R^N}m\left| \frac{w}{m}\right|^{\g'}dx<+\infty$. Then,
\begin{enumerate}
    \item[i)] $$m\in L^\beta(\R^N),\quad\forall \beta\in\left[1,\frac{N}{N-\g'}\right)\qquad(\forall\beta\in[1,+\infty), \,\,\,\text{if}\,\,\, \g'\ge N)$$
    and there exists a constant $C$ depending on $N$, $\beta$ and $\gamma'$ such that $$\|m\|_{L^\beta(\R^N)}\le C(E+M)$$
    \item[ii)] $$m\in W^{1,\ell}(\R^N),\quad\forall \ell<q$$
    and there exists a constant $C$ depending on $N$, $\ell$ and $\g'$ such that $$\|m\|_{W^{1,\ell}(\R^N)}\le C(E+M);$$
    \item[iii)] if $\g'>N$, we have also
    $$m\in C^{0,\theta}(\R^N),\quad \forall\theta\in (0, 1-N/\g')$$ 
    and there exists a constant $C$ depending on $N$, $\theta$ and $\g'$ such that
    $$\|m\|_{C^{0,\theta}}(\R^N)\le C(E+M).$$
\end{enumerate}
Moreover, we have
\begin{equation}\label{stima m_p 1}
\|m\|_{L^\frac{2N}{N+\al}(\R^N)}^\frac{2\g'}{N-\al}\le CM^{\frac{2\g'} {N-\al}-1}E
\end{equation}
where $C$ is a constant depending on $N,\g$ and $\al$;
\end{prop}

\proof
For the proof of statement $i)$ and $ii)$ see \cite[Proposition 2.5]{BerCes}. \textit{Proof of $iii)$.} From $ii)$ we have in particular that $m\in W^{1,\ell}(\R^N)$ for $N<\ell<\g'$, hence by Morrey's embedding
$$m\in C^{0,\theta}(\R^N),\quad\text{for}\,\,\,0<\theta<1 -\frac{N} {\g'}$$
and there exists a constant $C$, depending on $\theta$, $N$ and $\g'$, such that
$$\|m\|_{C^{0,\theta}(\R^N)}\le C(E+M).$$
Concerning estimate \eqref{stima m_p 1} refer to \cite[Proposition 2.4]{BerCes}.
\endproof

\noindent The following compactness result for sequences of couples $(m_n,w_n)\in \K_{\e,M}$ holds.

\begin{prop}\label{conv_mn}
Let us consider a sequence of couples $(m_n,w_n)\in\K_{\e,M}$ such that $E_n<C$ uniformly in $n$. Assume also that there exists a couple $(\bar{m},\bar{w})\in \K_{\e,M}$ such that $\bar{E}<+\infty$ and $m_n\to\bar{m}$ in $L^1(\R^N)$ as $n\to+\infty$. Then
$$m_n\to\bar{m}\,\,\,in\,\,\,L^s(\R^N),\quad\forall s\in\left[1,\frac{N} {N-\g'}\right)$$
(the previous holds $\forall s\in\left[1,+\infty\right)$ if $\g'\ge N$).
\end{prop}

\proof
If $\g'<N$, from Proposition \ref{kolm} we have that $m_n,\bar{m}\in L^\beta(\R^N)$ $\forall\beta<\frac{N}{N-\g'}$ and  $\|m\|_{L^\beta(\R^N)}\le C(E+M)$. 
Let us pick $1\le s<\frac{N}{N-\g'}$ and $s_1\in\left(s,\frac{N} {N-\g'}\right)$, by interpolation we get there exists $\theta\in(0,1)$ (depending on $s$ and $s_1$) such that 
$$\|\bar{m}-m_n\|_{L^s(\R^N)}\le\|\bar{m}-m_n\|^\theta_{L^1(\R^N)}\,\|\bar{m}-m_n\|^{1-\theta}_{L^{s_1}(\R^N)}.$$
We observe that 
$$\|\bar{m}-m_n\|_{L^{s_1}(\R^N)}\le\|\bar{m}\|_{L^{s_1}}+\|m_n\|_{L^{s_1}}\le C(\bar{E} +M)+C(E_n+M)\le C_1.$$
Hence  $\|\bar{m}-m_n\|_{L^{s_1}(\R^N)}$ is bounded, since $\|\bar{m}-m_n\|_ {L^1(\R^N)}\to0$ we can conclude. The same argument holds in the case when $\g'\ge N$.
\endproof

Finally, we will need some a priori regularity estimates for solutions to the Hamilton-Jacobi-Bellman equation
\begin{equation}\label{HJB}
-\Delta u+\frac{1}{\g}|\nabla u(x)|^\g+\la=V(x)-K_\al\ast m(x)\quad\mathrm{in} \,\,\R^N,
\end{equation}
where $\g>1$ is fixed. Assuming $m\in L^1(\R^N)$ fixed and such that $K_\al\ast m$ is H\"older continuous, we can define 
\begin{equation}\label{lambda bar}
\bar{\la}:=\sup\{\la\in\R\,|\,\eqref{HJB}\,\,\,\text{has a solution}\,\,\,u\in C^2(\R^N)\}.
\end{equation}
If $\bar{\la}$ does exist, then there exists $u\in C^2(\R^N)$ solving the HJB equation with such value $\bar{\la}$, moreover $u$ is coercive and its gradient has polynomial growth. We refer the reader to \cite[§2.3]{BerCes} and the references \cite{BM, Ci, Ic1} for more details.

\subsection{Properties of the Riesz potential.}

The Riesz potential $K_\al$ is well-defined as an operator on the whole space $L^r(\R^N)$ if and only if $r\in\left[1,\frac{N}{\al}\right)$. We state now the following well-known theorem (for which refer e.g. to \cite[Theorem 14.37]{WZ} and \cite[Theorem 4.3]{LiebLoss}).

\begin{theorem}[Hardy-Littlewood-Sobolev inequality] \label{HLS}
Let $0<\al<N$ and $1<r<\frac{N}{\al}$. Then for any $f\in L^r(\R^N)$ it holds
$$\|K_\al\ast f\|_{L^{\frac{Nr}{N-\al r}}(\R^N)}\le C\|f\|_{L^r(\R^N)}$$
where $C$ is a constant depending only on $N$, $\al$ and $r$.
Moreover, let $s,t>1$ such that $\frac{1}{s}-\frac{\al}{N}+\frac{1}{t}=1$ and assume that $f\in L^s(\R^N)$ and $g\in L^t(\R^N)$. Then, there exists a sharp constant $C(N,\al,s)$ (independent of $f$ and $g$) such that
\begin{equation}\label{disug hls}
\left|\int_{\R^N}\int_{\R^N}\frac{f(x)\,g(y)}{|x-y|^{N-\al}}dx\,dy\right|\le C\|f\|_{L^s(\R^N)}\|g\|_{L^t(\R^N)}.
\end{equation}
\end{theorem}

Taking advantage of the previous result, we are able to prove more precise integrability results for the Riesz term $K_\al\ast m$ and a compactness result. 

\begin{corollary}\label{K_al in Lp e compattezza}
Assume that $(m,w)\in\K_{\e,M}$ and $E<+\infty$. If $\g'\ge N$ or $\g'<N$ and $N-\g'\le\al<N$, then 
$$K_\al\ast m\,\in L^\beta(\R^N),\quad\forall\beta \in\left(\frac{N}{N-\al},+\infty\right)$$
and there exists a constant $C$ depending on $N$, $\al$, $\g'$ and $\beta$ such that 
$$\|K_\al\ast m\|_{L^\beta(\R^N)}\le C(E+M).$$
Moreover, if we consider a sequence of couples $(m_n,w_n)\in\K_{\e,M}$ such that $E_n<C$ uniformly in $n$ and we assume also that there exists a couple $(\bar{m},\bar{w})\in \K_{\e,M}$ such that $\bar{E}<+\infty$ and $m_n\to\bar{m}$ in $L^1(\R^N)$ as $n\to+\infty$. It holds
\begin{equation}\label{compactness K_al m m }
(K_\al\ast m_n)m_n\xrightarrow[n\to+\infty]{}(K_\al\ast \bar{m})\bar{m},\quad \text{in}\,\,\,L^1(\R^N).
\end{equation}
\end{corollary}

\proof
\textit{Case $\g'\ge N$.} From Proposition \ref{kolm} i) we have in particular that $m\in L^\beta(\R^N)$ $\forall\beta\in(1,\frac{N}{\al})$, hence by Theorem \ref{HLS} it follows claim $i)$. \textit{Case $\g'<N$.} From Proposition \ref{kolm} i) it holds that $m\in L^\beta(\R^N)$ $\forall\beta< \frac{N}{N-\g'}$. In the case when $N-\g'\le\al<N$, we have that $m\in L^\beta(\R^N)$ $\forall\beta\in(1,\frac{N}{\al})$ and we can conclude as before. \textit{Proof of \eqref{compactness K_al m m }}. Let us consider $\bar{r}:=\frac{Nr}{N-\al r}$ and $(\bar{r})'$ its conjugate exponent, namely $(\bar{r})'=\frac{Nr}{Nr-N+\al r}$. If $\g'<N$, from Proposition \ref{conv_mn}, we observe that in order to have $m_n\to\bar{m}$ in $(L^r\cap L^{(\bar{r})'})(\R^N)$ for a certain $r\in\left(1,\frac{N}{N-\g'}\right)$, it is sufficient to require that $(\bar{r})'<\frac{N}{N-\g'}$, that is 
$$\frac{N}{\al+\g'}<r<\frac{N}{N-\g'}$$
and hence 
$$\al>N-2\g'.$$
In particular $m_n\to\bar{m}$ in $L^r(\R^N)$, so from Theorem \ref{HLS} it follows that $$K_\al\ast m_n\to K_\al\ast\bar{m},\quad\text{in}\,\,\,L^{\bar{r}}(\R^N)$$
and since $m_n\to\bar{m}$ in $L^{(\bar{r})'}(\R^N)$
$$(K_\al\ast m_n)m_n\to(K_\al\ast\bar{m})\bar{m},\quad\text{in}\,\,\,L^1(\R^N).$$
The case $\g'\ge N$ is analogous.
\endproof

Regarding the H\"older continuity and the $L^\infty$-norm of the Riesz potential, we remind here the following results (see e.g. \cite[Theorem 2.8]{BerCes}).

\begin{theorem}\label{holderRiesz}
Let $1<r<+\infty$ and $0<\al<N$ be such that $0<\al-\frac{N}{r}<1$. Then, for every $f\in L^1(\R^N)\cap L^r(\R^N)$ we have that 
$$K_\al\ast f\in C^{0,\al-\frac{N}{r}}(\R^N)$$
and there exists a constant $C$, depending on $r,\,\al$ and $N$, such that
$$\frac{\big|K_\al\ast f(x)-K_\al\ast f(y)\big|}{\|x-y\|^{\al-\frac{N}{r}}}\le C \|f\|_{L^r(\R^N)}.$$
\end{theorem}

\begin{theorem}\label{inftyRiesz}
Let $0<\al<N$, $1<r\le+\infty$ be such that $r>\frac{N}{\al}$ and $s\in[1,\frac{N} {\al})$. Then, for every $f\in L^s(\R^N)\cap L^r(\R^N)$ we have that 
\begin{equation}\label{K infty}
   \|K_\al\ast f\|_{L^\infty(\R^N)}\le C_1\|f\|_{L^r(\R^N)}+C_2\|f\|_{L^s(\R^N)} 
\end{equation}
where $C_1=C_1(N,\al,r)$ and $C_2=C_2(N,\al,s)$. 
\end{theorem}

Taking advantage of the integrability results in Proposition \ref{kolm}, we are able to prove H\"older continuity of the term $K_\al\ast m$, for couples $(m,w)\in\K_{\e,M}$ with finite kinetic energy $E$.

\begin{corollary}\label{cor_holder}
Assume that $(m,w)\in\K_{\e,M}$ and $E<+\infty$.
\begin{itemize}
    \item[i)] If $\g'\ge N$, then 
    $$K_\al\ast m\,\in C^{0,\theta}(\R^N),\quad\forall\theta \in(0,\min\{1,\al\}).$$
    \item[ii)] If $\g'<N$ and $\al>N-\g'$, then 
    $$K_\al\ast m\,\in C^{0,\theta}(\R^N),\quad \forall\theta \in(0,\min\{1,\al-(N-\g')\}\,).$$ 
\end{itemize} 
\end{corollary}

\proof
The thesis follows from Theorem \ref{holderRiesz} and Proposition \ref{kolm}. 
\endproof

We recall here a Brezis-Lieb-type lemma for the Riesz potential (refer to \cite[Theorem 1]{BL} for the classical Brezis-Lieb lemma). It will be a key tool in Section \ref{sez5}.

\begin{lemma}[Lemma 2.4 in \cite{MVs}]\label{brezis lieb}
Let $0<\al<N$, $p\in\left[1,\frac{2N}{N+\al}\right)$ and $(f_n)_{n\in\mathbb{N}}$ be a bounded sequence in $L^\frac{2Np}{N+\al}(\R^N)$. If $f_n\to f$ almost everywhere in $\R^N$ as $n\to+\infty$, then
$$\lim\limits_{n\to\infty}\int_{\R^N}(K_\al\ast|f_n|^p)|f_n|^p-\int_{\R^N}(K_\al\ast|f_n-f|^p)|f_n-f|^p=\int_{\R^N}(K_\al\ast|f|^p)|f|^p.$$
\end{lemma}

\subsection{Uniform $L^\infty$-bounds on $m$}

We recall here the following result, which provides uniform a priori $L^\infty$ bounds on $m$.

\begin{theorem}\label{teo4.1new}
Let $\al\in(N-\g',N)$ and $(s_k)_k,\,(t_k)_k$ be bounded positive real sequences. We consider a sequence of classical solutions $(u_k,m_k,\la_k)$ to the following MFG system
\begin{equation}\label{MFGk}
\begin{cases}
-\Delta u+\frac{1}{\g}|\nabla u|^{\g}+\lambda=s_k V(t_k x)-K_\al\ast m(x)\\
-\Delta m-\mathrm{div} \left(m\nabla u |\nabla u|^{\g-2} \right)=0\\
\int_{\R^N}m=M,\quad m\ge0
\end{cases}\text{in}\,\,\,\R^N,
\end{equation}
where the potential $V$ satisfies assumption \eqref{V} with constant $C_V,b$ independent of $k$. We assume that for every $k$, $m_k\in L^\infty(\R^N)$ and $u_k$ are bounded from below. Then, there exists a positive constant $C$ not depending on $k$ such that 
$$\|m_k\|_{L^\infty(\R^N)}\le C, \quad \forall k\in\mathbb{N}.$$
\end{theorem}

\proof
We follow the argument of the proof of \cite[Theorem 2.12]{BerCes} and \cite[Theorem 4.1]{CC}, but since we are in the mass-subcritical regime some assumptions we require in the former one can be weakened.\\
We assume by contradiction that 
$$\sup_{\R^N} m_k=L_k\to+\infty$$
and we define 
$$\de_k:=\begin{cases} 
L_k^{-\frac{1}{\al+\g'}} \quad\text{if }\,\,\g'\le N\\
L_k^{-\frac{1}{\g'}} \,\,\,\,\,\,\quad\text{if }\,\,\g'>N 
\end{cases}.$$
We rescale $(u_k,m_k,\la_k)$ as follows: 
$$v_k(x):=\de_k^{\frac{2-\g}{\g-1}}u_k(\de_k x)+1,\qquad n_k(x):=L_k^{-1}
m_k(\de_k x),\qquad \tilde{\la}_k:=\de_k^{\g'}\la_k.$$
In this way $0\le n_k(x)\le1$, $\sup n_k=1$ and $\int_{\R^N}n_k(x)dx=\de_k^a M\to0$ since $a=\al+\g'-N>0$ if $\g'\le N$ and $a=\g'-N>0$ if $\g'\ge N$. Moreover, since up to addition of constants we may assume $\inf u_k(x)=0$, we have $v_k(x)\geq 1$ for all $x\in\R^N$. We get that $(v_k, n_k,\tilde\la_k)$ is a solution to 
\begin{equation}\label{mfg vk mk}
\begin{cases}
-\Delta v_k+\frac{1}{\g}|\nabla v_k|^\g+\tilde{\la}_k=V_k(x)-K_\al\ast n_k(x)\\
-\Delta n_k-\mathrm{div}(n_k\nabla v_k|\nabla v_k|^{\g-2})=0
\end{cases}    
\end{equation}
where $$V_k(x):=\de_k^{\g'}s_kV(t_k\de_k x).$$
Notice that from \eqref{V}, denoting by $\s_k=\de_k^{\g'+b}s_k t_k^b$ we have
$$\de_k^{\g'}s_k C_V^{-1}(\max\{t_k\de_k|x|-C_V,0\})^b\le V_k(x)\le 
C_V(1+\s_k|x|^b),\quad \forall x\in\R^N$$
and by Theorem \ref{inftyRiesz} we get that
$$|K_\al\ast n_k(x)|\le C_{N,\al}\|n_k\|_\infty+\|n_k\|_1\le C+\de_k^a M\le 2C,\quad\text{uniformly in }k.$$
Moreover, by computing the HJB equation of \eqref{mfg vk mk} in a minimum point of $v_k$ we have that 
$$\tilde{\la}_k\ge-K_\al\ast n_k(\bar{x})\ge-C$$
while with the same computations as in \cite[Lemma 2.11]{BerCes} we obtain that
$\la_k\le C$ where $C$ depends on $\g, C_V, b,N$, this gives  $-C\le\tilde{\la}_k\le\delta_k^{\g'} C$ hence $|\tilde\la_k|\le C$ uniformly.

We can conclude following the proof of \cite[Theorem 4.1]{CC}. More in detail, if $x_k$ is an approximated maximum point of $n_k$ (that is $n_k(x_k)=1-\de$), then either $\s_k|x_k|^b\to+\infty$ up to subsequences or $\s_k|x_k|^b\le C$ for some $C>0.$
We suppose that the second possibility occurs, using a priori gradient estimates on $v_k$, we get that $n_k$ is uniformly (in $k$) H\"older continuous in the ball $B_1(x_k)$, which contradicts the fact that  $n_k\ge0$ and  $\|n_k\|_{L^1}\to 0$. 
Then $\s_k|x_k|^b\to+\infty$, in this case we may construct a Lyapunov function for the system and hence some integral estimates on $n_k$, this allows us to obtain a
uniform (in $k$) H\"older bound for $n_k$, which yields an absurd. This proves that $L_k\to +\infty$ is not possible.
\endproof


\section{Existence of ground states for $\e>0$}

In this section, we provide existence of classical solutions to the MFG system \eqref{1} using a minimization procedure. Notice that, even if this result partially covers the existence result obtained in \cite{BerCes}, the variational approach proves to be essential to obtain some suitable estimates that will be necessary in the vanishing viscosity setting and hence to prove concentration phenomena as $\e\to 0$. 

If $\g'<N$, condition
\begin{equation}\label{condmin}
N-\g'\,\le\,\,\al<N
\end{equation}
is \textit{necessary} for the problem 
$$\min\limits_{(m,w)\in\mathcal{K}_{\e,M}}\mathcal{E}(m,w)$$
to be well-posed. Indeed, let us consider $m=c e^{-|x|}$ such that $\int_{\R^N}m(x)dx=M$ and $w=\e\nabla m$, in this way $(m,w)\in\mathcal{K}_ {\e,M}$. 
For $\s>0$ define
$$m_\s(x):=\frac{m(\s^{-1}x)}{\s^N}\qquad\text{and}\qquad w_\s(x):= \frac{w(\s^{-1}x)}{\s^{N+1}},$$
we get that $(m_\s,w_\s)\in\mathcal{K}_{\e,M}$ and 
$$\E(m_\s,w_\s)=\frac{1}{\s^{\g'}}\left[\int_{\R^N}\frac{m}{\g'}\bigg|\frac{w}{m}\bigg|^{\g'}+\s^{\g'}Vm\,dX-\frac{\s^{\g'-N+\al}}{2}\int_{\R^N}m(X)(K_\al\ast m)(X)\,dX\right].$$
From \eqref{V} we have that $\int_{\R^N}\frac{m}{\g'}|\frac{w}{m}|^{\g'} +\s^{\g'}Vm\,dX\le C$, so if $\al<N-\g'$ and $\s\to0$ the Riesz term in the energy dominates and
$$\mathcal{E}(m_\s,w_\s)\to-\infty, \quad \text{as}\,\,\,\s\to 0.$$
Actually, condition \eqref{condmin} is also \textit{sufficient}, in fact we prove that if $N-\g'<\al<N$ the energy $\E$ is bounded from below; and in the case when $\al=N-\g'$, requiring in addition that the constraint mass $M$ is sufficiently small, the energy $\E$ is non-negative (see Section \ref{sez al=N-g'}). Hence, the minimum problem is well defined and by means of classical direct methods we are able to obtain minimizers. Notice that in the case when $\g'\ge N$, the above condition \eqref{condmin} reduces to $0<\al<N$.

In what follows we address the case $N-\g'<\al<N$. Without loss of generality, we may assume $\e\in(0,1]$ fixed. Let us define 
$$e_\e(M):=\inf\limits_{(m,w)\in\K_{\e,M}} \E(m,w).$$

\begin{lemma}\label{E fin}
Assume that $N-\g'<\al<N$ and let $(m,w)\in\K_{\e,M}$. Then, there exist $C_1=C_1(N,\g,\al,M)$, $C_2=C_2(N,\g,\al,M)$ and $K=K(M,C_V,b,N,\al,\g)$ positive constants such that
\begin{equation}\label{stima_basso_E}
-C_1\e^{-\frac{\g'(N-\al)}{\g'-N+\al}}\le e_\e(M) \le -C_2\e^{-\frac{\g'(N-\al)}{\g'-N+\al}}+K.
\end{equation}
\end{lemma}

\proof
Let us fix $\beta:=\frac{2N}{N+\al}$, since $1<\beta<1+\frac{\g'}{N}$ by \eqref{stima m_p 1}, \eqref{disug hls} and the fact that $V\ge0$, we get
\begin{equation}\label{E>=}
\E(m,w)\ge c_1 \e^{\g'}\|m\|_{L^\beta(\R^N)}^{\frac{2\g'}{N-\al}}-c_2\|m\|^2_{L^\beta(\R^N)}
\end{equation}
where $c_1$ is a constant depending on $N,\al,\g,M$ and $c_2$ is a constant which depends on $N$ and $\al$. Minimizing the RHS of \eqref{E>=}, we obtain that 
$$c_1 \e^{\g'}\|m\|_{L^\beta(\R^N)}^{\frac{2\g'}{N-\al}}-c_2\|m\|^2_{L^\beta(\R^N)}\ge(c_3-c_4)\e^{-\frac{\g'(N-\al)}{\g'-N+\al}}$$
hence, there exists a constant $C_1>0$ depending on $N,\g,\al,M$ such that
$$\E(m,w)\ge-C_1\e^{-\frac{\g'(N-\al)}{\g'-N+\al}}.$$
In order to prove the estimate from above it is enough to show that for a suitable couple $(\tilde{m},\tilde{w})\in\K_{\e,M}$ it holds
\begin{equation}\label{4.24}
\E(\tilde{m},\tilde{w})\le-C_2\e^{-\frac {\g'(N-\al)}{\g'-N+\al}}+K.
\end{equation}
Let us consider a smooth function $\p:[0,+\infty)\to\R$ defined as $\p(r)=e^{-r}$. We define 
$$\tilde{m}(x):=M\tau^N I_1\p(\tau|x|)$$
$$\tilde{w}(x):=\e\nabla\tilde{m}(x)$$
where $\tau$ is a positive constant to be fixed and $I_1^{-1}:=\int_{\R^N} e^{-|y|}dy$, obviously $(\tilde{m},\tilde{w})\in\K_{\e ,M}$. We get that
$$\int_{\R^N}\tilde{m}\left|\frac{\tilde{w}}{\tilde{m}}\right|^{\g'}dx=M\e^{\g'}\tau^{\g'},$$
$$\int_{\R^N}V(x)\,\tilde{m}\,dx\le MC_V+MC_V I_1 I_2\frac{1}{\tau^b}$$
and
$$\int_{\R^N}\int_{\R^N}\frac{\tilde{m}(x)\tilde{m}(y)}{|x-y|^{N-\al}}dx\,dy=M^2 I_1^2 I_3 \tau^{N-\al}$$
where $I_2:=\int_{\R^N}|y|^b\p(|y|)dy$ and $I_3:=\int_{\R^N}\int_{\R^N}\frac{\p(|x|) \p(|y|)}{|x-y|^{N-\al}}dx\,dy$. Now, coming back to the energy, we obtain
$$\E(\tilde{m},\tilde{w})\le M(\e\tau)^{\g'}+MC_V+M C_V I_1 I_2\frac{1}{\tau^b} -\frac{1}{2}M^2 I_1^2 I_3\,\tau^{N-\al}$$
finally taking $\tau=\frac{1}{A}\e^{-\frac{\g'}{\g'-N+\al}}$ we get
$$\E(\tilde{m},\tilde{w})\le\left(M\frac{1}{A^{\g'}}-\frac{1}{2}M^2 I_1^2 I_3\,\frac{1}{A^{N-\al}} \right) \e^{-\frac{\g'(N-\al)}{\g'-N+\al}}+MC_V+M C_V I_1 I_2\,\frac{1}{\tau^b}\le-C_2\e^{-\frac {\g'(N-\al)}{\g'-N+\al}}+K$$
choosing $A$ large enough.
\endproof

In particular, from Lemma \ref{E fin}, it follows that for $N-\g'<\al<N$, $e_\e(M)$ is finite. We have the following a priori bounds.

\begin{prop}\label{prop3.2}
Let $(m,w)\in\K_{\e,M}$ such that $e_\e(M)\ge\E(m,w)-\eta$ for some positive $\eta$. Then 
\begin{equation}\label{6.2}
\|m\|^2_{L^\frac{2N} {N+\al}(\R^N)}\le C\e^{-\frac{\g'(N-\al)} {\g'-N+\al}}
\end{equation}
\begin{equation}\label{6.2.1}
\int_{\R^N}m\left|\frac{w}{m}\right|^{\g'}dx\le C\e^{-\frac{\g' (N-\al)}{\g'-N+\al}}+K
\end{equation}
and
\begin{equation}\label{6.2.2}
\int_{\R^N}V(x)\,m\,dx\le C\e^{-\frac{\g'(N-\al)}{\g'-N+\al}}+K
\end{equation}
where $C$ and $K$ are positive constants depending on $M,N,\al, C_V, b,\g$ and $\eta$.
\end{prop}

\proof
Let us denote $\beta:=\frac{2N}{N+\al}$, if $(m,w)\in\K_{\e,M}$ and $e_\e(M)\ge\E (m,w)-\eta$ for some $\eta>0$, we have 
\begin{align}\label{aa}
c+\eta\ge e_\e(M)+\eta&\ge\E(m,w)\ge\frac{1}{\g'}\int_{\R^N} m\left|\frac{w}{m}\right|^{\g'}dx-\frac{1}{2}\int_{\R^N}\int_{\R^N}\frac{m(x)\,m(y)}{|x-y|^{N-\al}}dx\,dy\\ \notag
&\ge C_1\e^{\g'}M^{1-\frac{2\g'}{N-\al}}\|m\|_{L^\beta(\R^N)}^{\frac{2\g'}{N-\al}} -C_2\|m\|^2_{L^\beta(\R^N)},\quad 
\end{align}
where in the first inequality, we used that by Lemma \ref{E fin} there exists a positive constant $c$ depending on $M,C_V,b,\g,N, \al$ such that $e_\e(M)\le c$, while in the last inequality we exploit estimates \eqref{stima m_p 1} and \eqref{disug hls}. Since $\frac{\g'}{N-\al}>1$, choosing $C$ sufficiently large (not depending on $\e$) we have
$$C_1\e^{\g'}M^{1-\frac{2\g'}{N-\al}}\left(\e^{-\frac{\g'(N-\al)}{\g'-N+\al}}C\right)^\frac{\g'}{N-\al}-C_2\left(\e^{-\frac{\g'(N-\al)}{\g'-N+\al}}C\right)\ge c+\eta,$$
hence we must have 
$$\|m\|^2_{L^\frac{2N}{N+\al}(\R^N)}\le\e^{-\frac{\g'(N-\al)}{\g'-N+\al}}C.$$
From \eqref{aa} we get that 
\begin{align*}
\frac{1}{\g'}\int_{\R^N} m\left|\frac{w}{m}\right|^{\g'}dx&\le c+\eta+\frac{1}{2}\int_{\R^N}\int_{\R^N}\frac{m(x)\,m(y)}{|x-y|^{N-\al}}dx\,dy\\
&\le c+\eta+C \|m\|^2_{L^\frac{2N}{N+\al}(\R^N)}\le C\e^{-\frac{\g'(N-\al)  }{\g'-N+\al}}+K
\end{align*}
which proves \eqref{6.2.1}. Finally, since $\frac{1}{\g'}\int\limits_{\R^N}m \left|\frac{w}{m}\right|^{\g'}dx+\int\limits_{\R^N}V(x)\,m\,dx=\E(m,w) +\frac{1}{2}\int_{\R^N}\int_{\R^N}\frac{m(x)m(y)}{|x-y|^{N-\al}}$, using \eqref{disug hls} and \eqref{6.2}, we obtain
\begin{align}\notag
\int_{\R^N}V(x)m\,dx\le\frac{1}{\g'}\int\limits_{\R^N}m\left|\frac{w}{m}\right|^{\g'}dx+\int\limits_{\R^N} V(x)\,m\,dx&\le\E(m,w)+C\|m\|^2_{L^\frac{2N}{N+\al}(\R^N)}\\ \label{68new}
&\le e_{\e}(M)+\eta+C\e^{-\frac{\g'(N-\al)}{\g'-N+\al}}
\end{align}
which gives estimate \eqref{6.2.2}.
\endproof

By means of classical direct methods, we prove that for every $\e,M>0$ there exists a minimizer $(m_\e,w_\e)\in\K_{\e,M}$ of the energy $\E$. 

\begin{prop}\label{prop3.3}
For every $\e>0$ and $M>0$, there exists a minimizer $(m_\e,w_\e)\in\K_{\e,M}$ of the energy $\E$, namely 
$$\E(m_\e,w_\e)=\inf\limits_{(m,w)\in\K_{\e,M}}\E(m,w).$$
For every minimizer $(m_\e,w_\e)$ of $\E$, we have that 
\begin{equation}\label{|m_e|}
\|m_\e\|_{L^\frac{2N}{N+\al}(\R^N)}\le C\e^{-\frac{\g'(N-\al)}{\g'-N+\al}}
\end{equation}
\begin{equation}\label{int m|w/m|}
\int_{\R^N}m_\e\bigg|\frac{w_\e}{m_\e}\bigg|^{\g'}dx\le C\e^{-\frac{\g'(N-\al)} {\g'-N+\al}}+K
\end{equation}
and
\begin{equation}\label{int Vm_e<}
\int_{\R^N}V(x)m_\e\,dx\le C\e^{-\frac{\g'(N-\al)}{\g'-N+\al}}+K,
\end{equation}
where $C>0$ and $K$ are two constants not depending on $\e$. Moreover, it holds
\begin{equation}\label{3.12}
m_\e(1+|x|)^b\in L^1(\R^N)\qquad\text{and}\qquad w_\e(1+|x|)^{\frac{b}{\g}}\in L^1(\R^N).
\end{equation}
\end{prop}

\proof 
Let us consider a minimizing sequence $(m_n, w_n)\in\K_ {\e,M}$, namely $\E(m_n,w_n)\to e_{\e}(M)$ as $n\to+\infty$. For $n$ sufficiently large $e_\e(M)\ge\E(m_n,w_n)-1$, so estimates \eqref{6.2}, \eqref{6.2.1} and \eqref{6.2.2} hold. By Proposition \ref{kolm}, using \eqref{6.2.1}, we get that
$$\|m_n\|_{W^{1,r}(\R^N)}\le C, \quad\,\,\forall r<q$$
where $C$ does not depend on $n$, hence by Sobolev compact embedding, up to subsequences $m_n\to m_\e$ in $L^s(K)$ for $1\le s<q^*:=\frac{qN}{N-q}$ and $K\subset\subset\R^N$. We observe that if $A\subset\R^N$ it holds
$$\int_{A}m_n(x)\,dx=\int_{\R^N}m_n(x)\chi_A(x)\,dx\le\|m_n\|_{\frac{2N}{N+\al}}\|\chi_A\|_\frac{2N}{N-\al};$$
hence using \eqref{6.2} we get that for every $\mu>0$ there exists $\delta_\mu>0$ such that 
$$\int_{A}m_n(x)\,dx\le\mu$$
for every $n$ and for any $A\subset\R^N$ such that $meas(A)<\delta_\mu$. Using estimate \eqref{6.2.2} and \eqref{V} we obtain that for $R>1$
$$C\ge\int_{\R^N}m_nV\,dx\ge\int_{B_R^c}m_nV\,dx\ge CR^b\int_{B_R^c}m_n(x)\,dx$$
namely, for every $\eta>0$ there exists $R>1$ such that $\int_{|x|>R}m_n(x)\,dx\le\eta$ for any $n$ (more precisely, for every $n$ greater than a certain value $n_0$). Thus, using the Vitali Convergence Theorem, up to extracting a subsequence, we have that
$$m_n\to m_\e \quad\text{in}\,\,\,L^1(\R^N)$$
and consequently 
$$\int_{\R^N}m_\e(x)\,dx=M.$$
Moreover, since by Sobolev embeddings $m_n$ are bounded in $L^s(\R^N)$ for every $s\in[1,q^*)$, we have also that $m_n\to\bar{m}$ in $L^\frac{2N}{N+\al}(\R^N)$. If we assume $\beta=\frac{2N}{N+\al}$, by H\"older inequality we get also 
$$\int_{\R^N}|w_n|^{\frac{\g'\beta}{\g'-1+\beta}}dx\le\left(\int_{\R^N}m_n^{1-\g'}|w_n|^{\g'}dx\right)^{\frac{\beta}{\g'-1+\beta}}\|m_n\|_{L^\beta(\R^N)}^\frac{\beta(\g'-1)}{\g'-1+\beta}$$
hence using \eqref{6.2} and \eqref{6.2.1}
$$w_n\rightharpoonup w_\e\quad\text{in}\,\,\,L^\frac{\g'\beta}{\g'-1+ \beta}(\R^N).$$
From \eqref{6.2}, \eqref{6.2.1} and \eqref{6.2.2} passing to the limit as $n\to+\infty$ and using Fatou's Lemma we obtain estimates \eqref{|m_e|}, \eqref{int m|w/m|} and \eqref{int Vm_e<}.\\
\indent We can infer that $(m_\e,w_\e)\in\K_{\e,M}$. Since the functional $\int_{\R^N}m|\frac{w}{m}|^{\g'}+V(x)m\,dx$ is sequentially lower semi-continuous with respect to the weak convergence and by Corollary \ref{K_al in Lp e compattezza} we have that $(K_\al\ast m_n)m_n\to(K_\al\ast m_\e)m_\e$ in $L^1(\R^N)$, then $(m_\e,w_\e)$ is a minimum of the energy $\E$.\\
\indent Finally the fact that $m_\e(1+|x|)^b\in L^1(\R^N)$ follows from \eqref{int Vm_e<} and \eqref{V}; whereas, by H\"older inequality 
$$\int_{\R^N}|w_\e|(1+|x|)^{b/\g}dx\le\left(\int_{\R^N}m_\e^{-\frac{\g'}{\g}}|w_\e|^{\g'}dx\right)^\frac{1}{\g'}\left(\int_{\R^N}m_\e(1+|x|)^b dx\right)^\frac{1}{\g}$$
since $\int_{\R^N}m_\e^{-\frac{\g'}{\g}}|w_\e|^{\g'}dx=\int_{\R^N}m_\e\big|\frac{w_\e}{m_\e}\big|^{\g'}dx$, using \eqref{int m|w/m|} and the fact that $m_\e(1+|x|)^b\in L^1$, we obtain that $w_\e(1+|x|)^ {b/\g}\in L^1(\R^N)$. 
\endproof

Once we have obtained minimizers $(m_\e,w_\e)\in\K_{\e,M}$ of the energy $\E$, we construct the associated solutions $(u_\e,m_\e,\la_\e)$ of the MFG system \eqref{1}.

\begin{prop}\label{prop3.4}
Let $N-\g'<\al<N$. Assume that the potential $V$ is locally H\"older continuous and satisfies \eqref{V}. Then, for every $\e,M>0$ there exists a classical solution $(u_\e,m_\e,\la_\e)$ to the MFG system \eqref{1} such that 
$$u_\e(x)\ge C_\e|x|^{1+\frac{b}{\g}}-C_\e^{-1}$$
$$|\nabla u_\e(x)|\le C_\e(1+|x|^{\frac{b}{\g}})$$
where $C_{\e}$ positive constant. Moreover, $u_\e$ is unique up to additive constants, $m_\e\in L^\infty(\R^N)$ and there exist $C_1,C_2,K$ positive constants not depending on $\e$ such that 
\begin{equation}\label{stima lambda 2}
-C_1\e^{-\frac{\g'(N-\al)}{\g'-N+\al}}\le\la_\e\le K-C_2\e^{-\frac{\g'(N-\al)} {\g'-N+\al}}.
\end{equation}
\end{prop}

\proof
Following the proof of \cite[Proposition 3.4]{CC}, let us consider the set of test functions
$$\mathcal{A}:=\left\{\psi\in C^2(\R^N)\,\,\bigg|\,\, \limsup\limits_{|x|\to+\infty}\frac{|\nabla \psi(x)|}{|x|^{b/\g}}<+\infty, \,\,\,\limsup\limits_{|x|\to+\infty}\frac{|\Delta \psi(x)|}{|x|^b}<+\infty\right \}.$$
From Proposition \ref{prop3.3} there exists at least one minimizer $(m_\e,w_\e)$ of the energy $\E$, and one can verify (using \eqref{3.12} and integrating by parts) that
\begin{equation}\label{3.18}
-\e\int_{\R^N}m_\e\Delta\psi\,dx=\int_{\R^N}w_\e\cdot\nabla\psi\,dx,\quad \forall\psi\in \mathcal{A}
\end{equation}
(see (3.18) in \cite{CC}, for details). Since every minimizer satisfies \eqref{3.12} and \eqref{3.18}, minimizing $\E$ on $\K_{\e,M}$ is equivalent to minimize $\E$ on the following constraint set
$$\K:=\left\{(m,w)\in(L^1\cap W^{1,r})\times L^{\frac{\g'\beta}{\g'+\beta-1}} (\R^N)\,\bigg|\,(m,w)\,\,\,\text{satisfies}\,\,\, \eqref{3.12},\,\,\eqref{3.18}, \,\, \int_{\R^N}m=M,\,\, m\ge0\right\}$$
where $r<q$. Now we prove that if $(m_\e,w_\e)$ is a minimizer of $\E$ on $\K$, then $(m_\e,w_\e)$ is also a minimizer of the functional 
\begin{equation}\label{Jtilde}
\tilde{\E}(m,w):=\frac{1}{\g'}\int_{\R^N}m\left|\frac{w}{m}\right|^{\g'}dx+\int_{\R^N}V(x)m\,dx-\int_{\R^N}\int_{\R^N}\frac{m(x)\,m_\e(y)}{|x-y|^{N-\al}}dx\,dy
\end{equation}
on $\K$. Define
$$\Phi(m,w):=\begin{cases}\frac{1}{\g'}\int_{\R^N}m\left|\frac{w}{m}\right|^{\gamma'}dx\quad\text{if}\,\,\,(m,w)\in\K\\
+\infty\qquad\qquad\quad\qquad\text{otherwise}
\end{cases}$$
and 
$$\Psi(m):=\int_{\R^N}V(x)m\,dx-\frac{1}{2}\int_{\R^N}\int_{\R^N}\frac{m(x)\,m(y)}{|x-y|^{N-\al}}dx\,dy,$$
we have
$$\E(m,w)=\Phi(m,w)+\Psi(m).$$
For any $(m,w)\in\K$ and $\lambda\in(0,1)$ we define $m_\la:=(1-\la)m_\e+\la m$ and $w_\la:=(1-\la)w_\e+\la w$, by minimality of $(m_\e,w_\e)$ it holds
\begin{equation}\label{fi}
\Phi(m_\la,w_\la)-\Phi(m_\e,w_\e)\ge\Psi(m_\e)-\Psi(m_\la)
\end{equation}
and by convexity of $\Phi$
\begin{equation}\label{fii}
\la\left(\Phi(m,w)-\Phi(m_\e,w_\e)\right)\ge\Phi(m_\la,w_\la)-\Phi(m_\e,w_\e).
\end{equation}
From \eqref{fi} and \eqref{fii} we obtain that
$$\la\left(\Phi(m,w)-\Phi(m_\e,w_\e)\right)\ge-\la\int_{\R^N}V(m-m_\e)dx+\la\int_{\R^N}\int_{\R^N}\frac{m_\e(y)\left(m(x)-m_\e(x)\right)}{|x-y|^{N-\al}}dy\,dx+o(\la),$$
dividing by $\la$ and letting $\la$ go to 0, we get 
$$-\int_{\R^N}V(m-m_\e)dx+\int_{\R^N}\int_{\R^N}\frac{m_\e(y)\left(m(x)-m_\e(x)\right)}{|x-y|^{N-\alpha}}dy\,dx\le \Phi(m,w)-\Phi(m_\e,w_\e)$$
for any $(m,w)\in\K$. Hence, the couple $(m_\e,w_\e)$ minimizes $\tilde{\E}$ on $\K$.\\ Let us consider the following functional
$$\mathcal{L}(m,w,\la,\psi):=\tilde{\E}(m,w)+\int_{\R^N}\e m\Delta\psi+w\nabla\psi-\la m\,dx+\la M.$$
One can easily verify that
$$\min\limits_{(m,w)\in\K}\tilde{\E}(m,w)=\min\limits_{(m,w)\in E}\sup\limits_ {(\la,\psi)\in\R\times\mathcal{A}}\mathcal{L}(m,w,\la,\psi)$$
where 
$$E:=\left\{(m,w)\in(L^1\cap W^{1,r})(\R^N)\times L^ \frac{\g'\beta}{\g'+\beta-1} (\R^N)\,\bigg|\,(m,w)\,\,\,\text{satisfies} \,\,\eqref{3.12}\right\}.$$
Proceeding as in \cite[Proposition 3.4]{CC}, by means of the Fan's minimax theorem (refer to Theorem 2.3.7 in \cite{BV}) and the Rockafellar interchange theorem (see \cite[Theorem 3A]{Rockafellar}), we get that
\begin{align*}
\min\limits_{(m,w)\in E}&\sup\limits_{(\la,\psi)\in\R\times\mathcal{A}}\mathcal{L} (m,w,\la,\psi)\\
&=\sup\limits_ {(\la,\psi)\in\R\times\mathcal{A}}\min\limits_{(m,w)\in E}\int\limits_{\R^N} \frac{m}{\g'}\left|\frac{w}{m}\right|^{\g'}+Vm-(K_\al\ast m_\e)m+\e m\Delta\psi+ w\nabla\psi-\la m\,dx+\la M\\
&=\sup\limits_{(\la,\psi)\in\R\times\mathcal{A}}\int\limits_{\R^N}\min\limits_{(m,w)\in\R_{\ge0}\times\R^N}\left(\frac{m}{\g'}\left|\frac{w}{m}\right|^{\g'}+Vm-(K_\al\ast m_\e)m+\e m\Delta\psi+ w\nabla\psi-\la m\right)dx+\la M.
\end{align*}
We observe now that
$$\min\limits_{(m,w)\in\R_{\ge0}\times\R^N}\left(\frac{m}{\g'}\left|\frac{w}{m}\right|^{\g'}+Vm-(K_\al\ast m_\e)m+\e m\Delta\psi+ w\nabla\psi-\la m\right)=0$$
if $\e\Delta\psi-\frac{1}{\g}|\nabla\psi|^\g-\la+V-K_\al\ast m_\e\ge0$, while it is $-\infty$ otherwise. This proves that
$$\min\limits_{(m,w)\in\K}\tilde{\E}(m,w)=\sup\limits_{(\la,\psi)\in B}\la\,M$$
where 
$$B:=\left\{(\la,\psi)\in\R\times\mathcal{A}\,\bigg|\,-\e\Delta\psi+\frac{1}{\g}|\nabla\psi|^\g+\la\le V-K_\al\ast m_\e\,\,\,\text{on}\,\,\R^N \right\}.$$
From Corollary \ref{K_al in Lp e compattezza} and Corollary \ref{cor_holder} we have that $K_\al\ast m_\e\in L^\beta(\R^N)\cap C^{0,\theta}(\R^N)$, hence by \cite[Proposition 2.10]{BerCes} there exists a couple $(u_\e,\la_\e)\in C^2(\R^N)\times\R$ such that 
\begin{equation}\label{ubarsol}
-\e\Delta u_\e+\frac{1}{\g}|\nabla u_\e|^\g+\la_\e=V(x)-K_\al\ast m_\e(x),\quad \text{on}\,\,\,\R^N
\end{equation}
where
$$\la_\e=\bar{\la}_\e:=\sup\bigg\{\la\in\R\,\bigg|\, -\e\Delta u+\frac{1}{\g}|\nabla u|^\g+\la=V-K_\al\ast m_\e\,\,\,\text{has a subsolution}\,\,\,u_\e\in C^2(\R^N)\bigg\}$$
and $u_\e$ is unique up to additive constants. Using also \cite[Proposition 2.9]{BerCes} we have
$$u_\e(x)\ge C_\e|x|^{\frac{b}{\g}+1}-C_\e^{-1},\qquad |\nabla u_\e(x)|\le C_\e(1+|x|)^{\frac{b}{\g}}.$$
Since
$$\e|\Delta u_\e|\le\frac{1}{\g}|\nabla u_\e|^\g+|\la_\e|+V(x)+K_\al\ast m_\e(x)\le C(1+|x|)^b,$$
it follows that $\limsup\limits_{|x|\to+\infty}\frac{|\Delta u_\e (x)|}{|x|^b}<+\infty$, hence $u_\e\in\mathcal{A}$. This proves that $\sup\limits_{(\la,\psi)\in B}\la\,M=\la_\e M$ and consequently
\begin{equation}\label{minlambda}
\la_\e M=\tilde{\E}(m_\e,w_\e)=\E(m_\e,w_\e)-\frac{1}{2}\int_{\R^N}\int_{\R^N}\frac{m_\e(x)m_\e(y)}{|x-y|^{N-\al}}dx\,dy.
\end{equation}
Now, by \eqref{minlambda}, \eqref{ubarsol} and \eqref{3.18} (since $u_\e\in\mathcal{A}$) we get that
\begin{align*}
0&=\int_{\R^N}\left(\frac{1}{\g'}\left|\frac{w_\e}{m_\e}\right|^{\g'}+V(x)-K_\al\ast m_\e-\la_\e\right)m_\e\,dx\\
&=\int_{\R^N}\left(\frac{1}{\g'}\left|\frac{w_\e}{m_\e}\right|^{\g'}-\e\Delta u_\e+\frac{1}{\g}|\nabla u_\e|^\g \right)m_\e\,dx=
\int_{\R^N}\left(\frac{1}{\g'}\left|\frac{w_\e}{m_\e}\right|^{\g'}+\frac{w_\e}{m_\e}\cdot\nabla u_\e+\frac{1}{\g}|\nabla u_\e|^\g\right)m_\e\,dx
\end{align*}
and we must have
$$\frac{w_\e}{m_\e}=-\nabla u_\e|\nabla u_\e|^{\g-2},\quad\text{on the set}\,\,\,\{m_\e>0\}.$$
We can conclude that $\e\Delta m_\e+\text{div}(m_\e\nabla u_\e|\nabla u_\e|^{\g-2})=0$ in weak sense. 

\textit{Proof of estimate \eqref{stima lambda 2}.} From \eqref{minlambda} we get that
\begin{equation}\label{lamb}
\la_\e=\frac{1}{M}\E(m_\e,w_\e)-\frac{1}{2M}\int_{\R^N}\int_{\R^N}\frac{m_\e(x)m_\e(y)}{|x-y|^{N-\al}}dx\,dy
\end{equation}
hence, by \eqref{stima_basso_E} and \eqref{disug hls} we get that
\begin{equation}\label{lambda>}
\la_\e\ge-c_1\e^{-\frac{\g'(N-\al)}{\g'-N+\al}}-c_2\|m_\e\|^2_{L^\frac{2N}{N+\al}(\R^N)} \ge-C_1\e^{-\frac{\g'(N-\al)}{\g'-N+\al}}
\end{equation}
using \eqref{6.2} in the last inequality. Moreover, from \eqref{lamb} we have
$$\la_\e\le\frac{1}{M}\E(m_\e,w_\e)=\frac{1}{M}\inf\limits_{(m,w)\in\K_{\e,M}} \E(m,w)$$
and using \eqref{stima_basso_E}, we conclude the proof of estimate \eqref{stima lambda 2}.

Finally, the function $m_\e\in L^\infty(\R^N)\cap W^{1,p}(\R^N)$ for every $p\ge1$, this follows proving that the function $u_\e^r$ for $r>1$ is a Lyapunov function for the stochastic process with drift $\nabla u_\e|\nabla u_\e|^{\g-2}$ and $m_\e$ density of the invariant measure associated to the process, then using some results in \cite{MPR} and Proposition \ref{prop2.4CC} (we refer the reader to the proof of \cite[Proposition 4.3 $iv)$]{BerCes} for more details). This proves that the triple $(u_\e,m_\e,\la_\e)$ is a classical solution to the MFG system \eqref{1}.
\endproof

\subsection{Existence of minima in the critical case $\al=N-\g'$}\label{sez al=N-g'}

In this subsection, assuming $\g'<N$, we study the critical case $\al=N-\g'$. In order to deal with the Riesz-interaction term, which is not H\"older continuous a priori in the critical case $\al=N-\g'$ (refer to Corollary \ref{cor_holder}), we first regularize the problem convolving the Riesz term with a standard symmetric mollifier. More precisely, we consider the following regularised version of the MFG system \eqref{1}
\begin{equation}\label{1 caso critico}
\begin{cases}
-\e\Delta u+\frac{1}{\g}|\nabla u|^\g+\la=V(x)-K_{\al}\ast\p_k\ast m(x)\\
-\e\Delta m-\mathrm{div}(m\nabla u|\nabla u|^{\g-2})=0\\
\int_{\R^N}m=M,\quad m\ge0
\end{cases}\text{in}\,\,\R^N
\end{equation}
where $(\p_k)_k$ is a sequence of standard symmetric mollifiers approximating the unit as $k\to+\infty$. We associate to \eqref{1 caso critico} the following energy
\begin{equation}\label{energy reg}
\mathcal{E}_k(m,w):=\begin{cases} \int\limits_{\R^N}\frac{m}{\g'} \left|\frac{w}{m}\right|^{\g'}+V(x)m\,dx-\frac{1}{2}\int\limits_{\R^N} m(x)K_\al\ast\p_k\ast m(x)dx\quad \text{if}\,\,(m,w)\in \mathcal{K}_{\e,M}\\
+\infty \hspace{8.1cm} \text{otherwise}\end{cases}
\end{equation}
and we prove that if the total mass $M$ is sufficiently small then, $\E_k$ is bounded from below.

\begin{lemma}\label{E fin 2}
Let $\al=N-\g'$ and $(m,w)\in\K_{\e,M}$. Then, there exists $M_0>0$ (depending on $N$, $\g$ and $\e$) such that for any $M\in(0,M_0]$
\begin{equation}\label{stima_basso_E 2}
\E_k(m,w)\ge0.
\end{equation}
Hence, $\inf_{(m,w)\in\K_{\e,M}}\E_k(m,w)$ is finite.
\end{lemma}

\proof 
Similarly to the proof of Lemma \ref{E fin}, let us fix $\beta:=\frac{2N}{2N-\g'}$, we get
\begin{equation}\label{E>= al=N-g'}
\E_k(m,w)\ge\left(C_1 \e^{\g'}\frac{1}{M}-C_2\right)\|m\|^2_{L^\beta(\R^N)}
\end{equation}
where $C_1$ and $C_2$ are constants depending on $N$ and $\g$. If $C_1 \frac{\e^{\g'}}{M}-C_2>0$, that is $M\le M_0$ where $M_0:=\frac{C_1}{C_2}\e^{\g'}$, we have that
$$\E_k(m,w)\ge0.$$
\endproof

\noindent As before, by classical direct methods we prove that for every $\e>0$ and $M\in(0,M_0]$ there exists a global minimizer $(m_k,w_k)\in\K_{\e,M}$ of the regularised energy $\E_k$, this allows us to construct the corresponding associated solutions $(u_k,m_k,\la_k)$ of regularised problem. With the same arguments of \cite[§4.2]{BerCes}, since we have a uniform $L^\infty$-bound on $m_k$, we can finally pass to the limit as $k\to+\infty$ in the MFG system and obtain a solution to the initial problem \eqref{1} for $\al=N-\g'$.


\section{Asymptotic analysis of solutions}\label{sez5}

In this section, assuming $\al\in(N-\g',N)$, we want to study the behavior of agents when the Brownian noise vanishes. To this end, we analyze the asymptotic behavior of a solution $(m_\e,u_\e,\la_\e)$ to the MFG system \eqref{1} as $\e\to0$. 

\subsection{The rescaled problem and some a priori estimates}\label{pb riscalato}

For $\e>0$, let us define a suitable rescaling for $m$, $u$ and $\la$, which preserves the mass of $m$:
$$\tilde{m}(y):=\e^{\frac{N\g'}{\g'-N+\al}}m\left(\e^{\frac{\g'}{\g'-N+\al}}y \right),$$
$$\tilde{u}(y):=\e^{\frac{\g'(N-\al)-\g'-N+\al}{\g'-N+\al}}u\left(\e^{\frac{\g'}{\g'-N+\al}}y\right),$$
$$\tilde{\la}:=\e^{\frac{(N-\al)\g'}{\g'-N+\al}}\la$$
and a rescaled potential
\begin{equation}\label{def V_e}
V_\e(y):=\e^{\frac{(N-\al)\g'}{\g'-N+\al}}V\left(\e^{\frac{\g'}{\g'-N+\al}}y\right).
\end{equation}
Notice that $V_\e$ vanishes locally as $\e\to0$, hence, passing to the limit, we can not take advantage of the coercivity of $V$ in order to prove that there is no loss of mass (compare with the proof of Proposition \ref{prop3.3}) indeed we will use a concentration-compactness argument. From assumptions \eqref{V} on the potential $V$, we get the corresponding assumptions on $V_\e$
\begin{equation}\label{V_e}
C_V^{-1}\e^{\frac{(N-\al)\g'}{\g'-N+\al}}\left(\max\left\{\e^{\frac{\g'}{\g'-N+\al}}|y|-C_V,0\right\}\right)^b\le V_\e(y)\le C_V\e^{\frac{(N-\al) \g'}{\g'-N+\al}}\left(1+\e^{\frac{\g'}{\g'-N+\al}}|y|\right)^b.
\end{equation}
The rescaled Riesz-type interaction term is
$$K_\al\ast\tilde{m}(y)=\e^{\frac{(N-\al)\g'}{\g'-N+\al}}K_\al\ast m\left(\e^{\frac {\g'}{\g'-N+\al}}y\right).$$
Hence, if the triple $(u_\e,m_\e,\la_\e)$ is a classical solution to the MFG system \eqref{1} (from Proposition \ref{prop3.4} there exists at least one solution to \eqref{1}), one can verify that 
\begin{align*}
-\Delta&\tilde{u}_\e(y)+\frac{1}{\g}|\nabla\tilde{u}_\e(y)|^\g+\tilde{\la}_\e\\
&=-\e^{\frac{\g'(N-\al)}{\g'-N+\al}+1}\Delta u_\e(\e^\frac{\g'}{\g'-N+\al}y)+
\e^{\frac{\g'(N-\al)}{\g'-N+\al}}\frac{1}{\g}\big|\nabla u_\e\big(\e^\frac{\g'}{\g'-N+\al}y\big)\big|^\g+\e^{\frac{\g'(N-\al)}{\g'-N+\al}}\la_\e\\
&=\e^\frac{\g'(N-\al)}{\g'-N+\al}\left[-\e\Delta u_\e+\frac{1}{\g}|\nabla u_\e|^\g+\la_\e\right]=\e^\frac{\g'(N-\al)}{\g'-N+\al}\left[V(\e^\frac{\g'}{\g'-N+\al}y)-K_\al\ast m(\e^\frac{\g'}{\g'-N+\al}y)\right]\\
&=V_\e(y)-K_\al\ast \tilde{m}(y)
\end{align*}
and also
\begin{align*}
-\Delta&\tilde{m}_\e(y)-\mathrm{div}\big(\tilde{m}_\e(y)\nabla\tilde{u}_\e(y)|\nabla \tilde{u}_\e(y)|^{\g-2}\big)\\
&=-\e^\frac{N\g'+2\g'}{\g'-N+\al}\Delta m_\e(\e^\frac{\g'}{\g'-N+\al}y)-\mathrm{div}\left(\e^\frac{N\g'+N-\al}{\g'-N+\al}m_\e\nabla u_\e|\nabla u_\e|^{\g-2}(\e^\frac{\g'}{\g'-N+\al}y)\right)\\
&=\e^\frac{N\g'+N-\al+\g'}{\g'-N+\al}[-\e\Delta m_\e-\mathrm{div}(m_\e\nabla u_\e|\nabla u_\e|^{\g-2})]=0.
\end{align*}
Hence
\begin{equation}\label{5.13}
\begin{cases}
-\Delta\tilde{u}_\e+\frac{1}{\g}|\nabla\tilde{u}_\e|^\g+\tilde{\la}_\e=V_\e(y)-K_\al\ast\tilde{m}_\e(y)\\
-\Delta\tilde{m}_\e-\mathrm{div}(\tilde{m}_\e\nabla \tilde{u}_\e|\nabla \tilde{u}_\e|^{\g-2})=0\\
\int_{\R^N}\tilde{m}_\e=M
\end{cases}\quad\text{in}\,\,\,\R^N.
\end{equation}
In order to prove that there is no loss of mass when passing to the limit as $\e\to0$, we translate the reference system at a point around which the mass $\tilde{m}_\e$ remains positive. By Proposition \ref{prop3.4}, $u_\e\in C^2(\R^N)$ is coercive, hence there exists a point $y_\e\in\R^N$ such that
$$\tilde{u}_\e(y_\e)=\min\limits_{\R^N}\tilde{u}_\e(y).$$
Let us define 
$$\bar{u}_\e(y):=\tilde{u}_\e(y+y_\e)-\tilde{u}_\e(y_\e)$$
$$\bar{m}_\e(y):=\tilde{m}_\e(y+y_\e)$$
in this way we have $\bar{u}_\e(0)=0=\min_{\R^N}\bar{u}_\e$. One can immediately verify that $(\bar{m}_\e,\bar{u}_\e,\tilde{\la}_\e)$ is a classical solution to 
\begin{equation}\label{5.19}
\begin{cases}
-\Delta\bar{u}_\e+\frac{1}{\g}|\nabla\bar{u}_\e|^\g+\tilde{\la}_\e=V_\e(y+y_\e)-K_\al\ast\bar{m}_\e(y)\\
-\Delta\bar{m}_\e-\mathrm{div}(\bar{m}_\e\nabla\bar{u}_\e|\nabla\bar{u}_\e|^{\g-2})=0\\\int_{\R^N}\bar{m}_\e=M
\end{cases}.
\end{equation}

We define also the rescaled energy
$$\E_\e(m,w):=\int_{\R^N}\frac{m}{\g'}\left|\frac{w}{m}\right|^{\g'}+V_\e(y+y_\e)m\,dy-\frac{1}{2}\int_{\R^N}\int_{\R^N}\frac{m(x)m(y)}{|x-y|^{N-\al}}dx\,dy.$$
Notice that $\E_\e(\bar{m}_\e,\bar{w}_\e)=\e^\frac{(N-\al)\g'}{\g'-N+\al}\E(m_\e,w_\e)$, hence if $(m_\e,w_\e)\in\K_{\e,M}$ is a minimizer of $\E$, then $(\bar{m}_\e,\bar{w}_\e)$ is a minimizer of $\E_\e$ on $\K_{1,M}$ (where $\bar{w}_\e=-\bar{m}_\e\nabla\bar{u}_\e|\nabla\bar{u}_\e|^{\g-2}$). We will denote 
$$\tilde{e}_\e(M):=\min\limits_{(m,w)\in\K_{1,M}}\E_\e(m,w).$$
From \eqref{stima_basso_E} by rescaling, we get
\begin{equation}\label{5.29}
-C_1\le\tilde{e}_\e(M)\le -C_2+K\e^{\frac{(N-\al)\g'}{\g'-N+\al}}
\end{equation}
where $C_1, C_2, K$ are positive constants not depending on $\e$.\\

First of all, we prove the following a priori estimates. 

\begin{lemma}
Let $(\bar{u}_\e,\bar{m}_\e,\tilde{\la}_\e)$ be a classical solution to \eqref{5.19}.
Then there exist $C_1,C_2,C$ positive constants not depending on $\e$ such that
\begin{equation}\label{5.14}
     -C_1\le\tilde{\la}_\e\le-C_2
\end{equation}

\begin{equation}\label{5.15}
     \int_{\R^N}\bar{m}_\e|\nabla\bar{u}_\e|^{\g}dx\le C
\end{equation}

\begin{equation}\label{5.16}
    \|\bar{m}_\e\|_{L^\infty(\R^N)}\le C
\end{equation}

\begin{equation}\label{5.20a}
    \e^\frac{\g'(N-\al+b)}{\g'-N+\al}|y_\e|^b\le C   
\end{equation}

\begin{equation}\label{5.20b}
    0\le V_\e(y+y_\e)\le C\left(1+\e^\frac{\g'(N-\al+b)}{\g'-N+\al}|y|^b\right)
\end{equation}

\begin{equation}\label{5.21}
|\nabla\bar{u}_\e(y)|\le C(1+|y|)^\frac{b}{\g}\qquad\text{and}\qquad \bar{u}_\e(y) \ge C|y|^{1+\frac{b}{\g}}-C^{-1}.
\end{equation}

\noindent Moreover, for $R$ sufficiently large we have
\begin{equation}\label{5.22}
\int_{B_R(0)}\bar{m}_\e(y)dy\ge C.
\end{equation}

\end{lemma}

\proof 
Estimates \eqref{5.14} and \eqref{5.15} follow, by rescaling, from \eqref{stima lambda 2} and \eqref{int m|w/m|} respectively. From Proposition \ref{prop3.4} we have that for every $\e$, $\bar{u}_\e$ are bounded from below and $\bar{m}_\e\in L^\infty(\R^N)$, so by Theorem \ref{teo4.1new} with $s_\e=\e^\frac{(N-\al)\g'}{\g'-N+\al}$ and $t_\e=\e^\frac{\g'}{\g'-N+\al}$, we obtain the uniform $L^\infty$-bound \eqref{5.16}. Evaluating the first equation of \eqref{5.19} in $y=0$, we get
\begin{equation}\label{1 eq in 0}
\tilde{\la}_\e\ge V_\e(y_\e)-K_\al\ast\bar{m}_\e(0)   
\end{equation}
from estimates \eqref{5.14}, \eqref{5.16} and \eqref{V_e} we get \eqref{5.20a}; using it and \eqref{V_e} again, we obtain \eqref{5.20b}.\\
Since \eqref{5.14}, \eqref{5.16} and \eqref{5.20b} hold and $\bar{u}_\e$ is bounded from below, from \cite[Proposition 2.9]{BerCes} we get estimates \eqref{5.21} (which are uniform with respect to $\e$ since $\|K_\al\ast\bar{m}_\e\|_\infty\le C_{N,\al}\|\bar{m}_\e\|_\infty+M\le C$ uniformly in $\e$).\\
From \eqref{1 eq in 0}, using the fact that $V_\e\ge0$ and \eqref{5.14}, we get that there exists a positive constant $C$ not depending on $\e$ such that
$$K_\al\ast\bar{m}_\e(0)\ge C>0,$$
hence
$$C\le\int_{B_R}\frac{\bar{m}_\e(y)}{|y|^{N-\al}}dy+\int_{\R^N\setminus B_R}\frac{\bar {m}_\e(y)}{|y|^{N-\al}}dy\le\int_{B_R}\frac{\bar{m}_\e(y)}{|y|^{N-\al}}dy+\frac{M}{R^{N-\al}}.$$
This implies that for $R>0$ sufficiently large 
$$\int_{B_R}\frac{\bar{m}_\e(y)}{|y|^{N-\al}}dy\ge C_1>0.$$
Moreover, if $r<R$ we have
$$C_1\le \int_{B_R\setminus B_r}\frac{\bar{m}_\e(y)}{|y|^{N-\al}}dy+\int_{B_r} \frac{\bar{m}_\e(y)}{|y|^{N-\al}}dy\le \frac{1}{r^{N-\al}}\int_{B_R\setminus B_r}\bar{m}_\e(y)dy+\|\bar{m}_\e\|_\infty\int_{B_r}\frac{dy}{|y|^{N-\al}}.$$
Keeping in mind that $\int_{B_r}\frac{dy}{|y|^{N-\al}}=cr^\al$ and \eqref{5.16}, we can infer that choosing $r$ sufficiently small
$$\frac{1}{r^{N-\al}}\int_{B_R\setminus B_r}\bar{m}_\e(y)dy\ge C_2>0$$
and consequently
$$\int_{B_R}\bar{m}_\e(y)dy\ge C_3>0$$
\endproof


\subsection{Convergence of solutions}

At this stage we are able to prove a convergence result, which however, do not ensure conservation of mass in the limit. 

\begin{prop}\label{prop5.3}
If $(\bar{u}_\e,\bar{m}_\e,\tilde{\lambda}_\e)$ is a classical solution to \eqref{5.19}, then as $\e\to0$ up to extracting a subsequence we have that 
$$\tilde{\la}_\e\to\bar{\la}$$
and
$$\bar{u}_\e\to\bar{u},\qquad\nabla\bar{u}_\e\to\nabla\bar{u},\qquad\bar{m}_\e\to\bar{m},\qquad\text{locally uniformly}.$$
The triple $(\bar{u},\bar{m},\bar{\la})$ is a classical solution to 
\begin{equation}\label{5.25}
\begin{cases}
-\Delta\bar{u}+\frac{1}{\g}|\nabla\bar{u}|^\g+\bar{\la}=g(x)-K_\al\ast\bar{m}(x)\\
-\Delta\bar{m}-\mathrm{div}(\bar{m}\nabla\bar{u}|\nabla\bar{u}|^{\g-2})=0\\
\int_{\R^N}\bar{m}\,dx=a
\end{cases}
\end{equation}
where $g$ is a continuous function such that, up to subsequence, $V_\e(x+y_\e) \xrightarrow[\e\to0]{}g(x)$ locally uniformly and $a\in(0,M]$. Moreover, there exist $C_1,C_2,c_1,c_2$ positive constants such that
\begin{equation}\label{5.26}
\bar{u}(y)\ge C_1|y|-{C_1}^{-1} \qquad\text{and}\qquad |\nabla\bar{u}|\le C_2
\end{equation}
and finally
\begin{equation}\label{barm<ce^-cx}
\bar{m}(x)\le c_1e^{-c_2|x|}, \quad\text{on}\,\,\,\R^N.
\end{equation}
\end{prop}

\proof 
By means of the previous uniform estimates and the fact $\bar{u}_\e$ is a classical solution to the HJB equation, using \cite[Theorem 8.32]{GT}, we get that for any compact set $K$ in $\R^N$ and for any $\theta\in(0,1]$ 
$$\|\bar{u}_\e\|_{C^{1,\theta}(K)}\le C\qquad\text{locally uniformly in }\e.$$
Using \eqref{5.16} and \eqref{5.21}, by Proposition \ref{prop2.4CC} and Sobolev embedding, we get that for every $\mu\in(0,1)$
$$\|\bar{m}_\e\|_{C^{0,\mu}(K)}\le C\qquad\text{locally uniformly in }\e.$$
Hence, up to subsequences, we have that as $\e\to0$
$$\bar{u}_\e\to\bar{u},\qquad\text{locally uniformly in $C^1$ on compact sets}$$
and 
$$\bar{m}_\e\to\bar{m} \qquad\text{locally uniformly on compact sets}$$
and weakly in $W^{1,p}(B_R)$ $\forall p>1$ and $R>0$. Moreover, from \eqref{5.14} it follows that up to extracting a subsequence $\tilde{\la}_\e\to\bar{\la}$. From \eqref{5.16} and the fact that $\|\bar{m}_\e\|_{L^1(\R^N)}=M$, by interpolation we get that for every $p\in(1,+\infty)$ it holds $\|\bar{m}_\e\|_{L^p(\R^N)}\le C$ uniformly in $\e$. 
Using Theorem \ref{holderRiesz} and Theorem \ref{inftyRiesz} we have that $\|K_\al\ast\bar{m}_\e\|_{C^{0,\al-\frac{N}{r}}(\R^N)}\le C$  uniformly in $\e$ and up to extracting a subsequence  
$$K_\al\ast\bar{m}_\e\xrightarrow[\e\to0]{} K_\al\ast\bar{m},\quad\text{locally uniformly}.$$
Moreover, from estimate \eqref{5.20b} we have that up to subsequences $V_\e(x+y_\e) \to g(x)$ locally uniformly, where $g$ is a continuous function such that $0\le g(x)\le C$. Finally, from \eqref{5.22} it follows that $\int_{\R^N}\bar{m}(x)\,dx=a\in(0,M]$.
By stability with respect to uniform convergence, $\bar{u}$ solves in the viscosity sense 
$$-\Delta u+\frac{1}{\g}|\nabla u|^\g+\bar{\la}=g(x)-K_\al\ast\bar{m}(x)$$
and using the regularity of the HJB equation we get that $\bar{u}\in C^2$. Finally, by the strong convergence $\nabla u_\e\to\nabla\bar{u}$, we get that $\bar{m}$ solves
$$\Delta m-\mathrm{div}(m\nabla\bar{u}|\nabla\bar{u}|^{\g-2})=0$$ and with the same procedure as \cite[Proposition 4.3 $iv)$]{BerCes} we get that $\bar{m}\in W^{1,p}(\R^N)$ for every $p\in(1,+\infty)$. Hence, $(\bar{u},\bar{m},\bar{\la})$ is a classical solution to \eqref{5.25}.\\
In order to prove \eqref{5.26} we use \cite[Proposition 2.9]{BerCes}. Notice that, if $f$ is a non-negative H\"older continuous function such that  $\int_{\R^N}f^\beta\, dx<+\infty$ for a certain $\beta>1$, then $f(x)\to0$ as \blu{$|x|\to+\infty$} (see \cite[Lemma 2.2]{CC} which is stated in the case $\beta=1$ but it can be easily generalised to $\beta>1$). Since $K_\al\ast\bar{m}\in C^{0,\theta}(\R^N)\cap L^\beta(\R^N)$ and it is non-negative, we get that 
$$K_\al\ast\bar{m}(x)\to0, \quad \text{as}\,\,\,|x|\to+\infty$$
and hence
$$\liminf\limits_{|x|\to+\infty}\left(g(x)-K_\al\ast\bar{m}(x)-\bar{\la}\right)\ge-\bar{\la}>0.$$
From \cite[Proposition 2.9]{BerCes} we get \eqref{5.26}. 

Since we can choose $k>0$ such that the function $\p(x):=e^{k\bar{u}(x)}$ is a Lyapunov function for the process, from \cite[Proposition 2.6]{MPR} we get that 
$$e^{k\bar{u}}\in L^1(\bar{m})$$
and finally from \cite[Theorem 6.1]{MPR} it follows \eqref{barm<ce^-cx}.
\endproof


\subsection{No loss of mass when passing to the limit}

First, we prove that the energy functional $\E_\e(m,w)$ holds a sort of sub-additive property. Then, we assume by contradiction to have loss of mass, namely that $\int_{\R^N}\bar{m}\,dx=a\in(0,M)$, by means of a concentration-compactness argument we prove that this leads to an absurd. Hence $\bar{m}$ has still $L^1$-norm equal to $M$.

\begin{lemma}\label{lemma5.5}
For all $a\in(0,M)$, there exists a positive constant $C$ depending on $a$, $M$ and the other constants of the problem (but not on $\e$) such that 
\begin{equation}
\tilde{e}_\e(M)<\tilde{e}_\e(a)+\tilde{e}_\e(M-a)-C.
\end{equation}
\end{lemma}

\proof
Let us assume that $a\ge\frac{M}{2}$ and fix $c>1$ and $B>0$. If $(m,w)\in\K_{1,B}$ we get
\begin{align}\notag
\tilde{e}_\e(cB)&\le\E_\e(cm,cw)=\int_{\R^N}\frac{cm}{\g'}\left|\frac{w}{m}\right|^{\g'}+cV_\e(x+y_\e)m\,dx-\frac{c^2}{2}\int_{\R^N}\int_{\R^N}\frac{m(x)\,m(y)}{|x-y|^{N-\al}}dx\,dy\\\label{5.31}
&=c\E_\e(m,w)-\frac{c(c-1)}{2}\int_{\R^N}\int_{\R^N}\frac{m(x)\,m(y)}{|x-y|^{N-\al}}dx\,dy
\end{align}
If $(m,w)\in\K_{1,B}$ is a minimizer of $\E_\e$, we have
$$-C_2(B)+K\e^\frac{(N-\al)\g'}{\g'-N+\al}\ge\tilde{e}_\e(B)\ge-\frac{1}{2}\int_{\R^N}\int_{\R^N}\frac{m(x)\,m(y)}{|x-y|^{N-\al}}dx\,dy$$
notice that the constant $C_2$ is the one that appears in \eqref{5.29} and depends on $B$ and on the others variables of the problem. Taking $\e$ sufficiently small, we obtain
\begin{equation}\label{5.31b}
\frac{1}{2}\int_{\R^N}\int_{\R^N}\frac{m(x)\,m(y)}{|x-y|^{N-\al}}dx\,dy\ge\frac{C_2(B)}{2}>0.
\end{equation}
and using \eqref{5.31b} in \eqref{5.31} we get
\begin{equation}\label{5.32}
\tilde{e}_\e(cB)<c\tilde{e}_\e(B)-c(c-1)\frac{C_2(B)}{2}.
\end{equation}
Taking $B=a$ and $c=M/a$ in \eqref{5.32} we get
$$\tilde{e}_\e(M)<\frac{M}{a}\tilde{e}_\e(a)-\frac{M}{a}\left(\frac{M}{a}-1\right)\frac{C_2(a)}{2}=\tilde{e}_\e(a)+\frac{M-a}{a}\tilde{e}_\e(a)-\frac{M}{a}\left(\frac{M}{a}-1\right)\frac{C_2(a)}{2}$$
if $a=M/2$ we have done, whereas if $a>M/2$ we take $B=M-a$ and $c=\frac{a}{M-a}$ in \eqref{5.32} and multiplying by $\frac{M-a}{a}$, we get
$$\frac{M-a}{a}\,\tilde{e}_\e(a)<\tilde{e}_\e(M-a)-\left(\frac{a}{M-a}-1\right)\frac{C_2(M-a)}{2}\le\tilde{e}_\e(M-a).$$
From the previous we can conclude that
$$\tilde{e}_\e(M)<\tilde{e}_\e(a)+\tilde{e}_\e(M-a)-\frac{M}{a}\left(\frac{M}{a}-1\right)\frac{C_2(a)}{2}.$$
In the case when $a<M/2$ we replace $a$ with $M-a$.
\endproof

From estimate \eqref{barm<ce^-cx}, it follows that there exists a positive constant $\bar{c}$ such that $\bar{m}\le \bar{c}e^{-|x|}$. For $R>0$ (which will be fixed later) let us define
$$\nu_R(x):=\begin{cases}\bar{c}e^{-R}\quad\,\,\text{if}\,\,\,|x|\le R\\
\bar{c}e^{-|x|}\quad\text{if}\,\,\,|x|>R \end{cases}.$$
We have the following splitting of the energy $\E_\e$.

\begin{lemma}\label{splitting >}
Let $(\bar{m}_\e,\bar{w}_\e)$ be a minimizer of $\E_\e$, $\bar{m}$ and $\bar{u}$ obtained from Proposition \ref{prop5.3} and $\bar{w}_\e\to\bar{w}=-\bar{m}\nabla \bar{u}|\nabla\bar{u}|^{\g-2}$ locally uniformly. If $\int_{\R^N}\bar{m}dx=a\in (0,M)$, then 
\begin{equation}\label{E_e>=}
\E_\e(\bar{m}_\e,\bar{w}_\e)\ge\E_\e(\bar{m},\bar{w})+ \E_\e(\bar{m}_\e-\bar{m}+2\nu_R,\bar{w}_\e-\bar{w}+2\nabla\nu_R)+o_\e(1)-C R^{b+N}e^{-R}
\end{equation}
as $\e\to0$.
\end{lemma}

\proof
Following the arguments of the proof of \cite[Theorem 5.6]{CC}, we recall some facts that we will need. By definition, $\bar{m}(x)\le\nu_R(x)$ for $|x|>R$ and
\begin{equation}\label{5.35}
\int_{\R^N}\nu_R(x)dx=\omega_NR^N\bar{c}e^{-R}+\int_{\R^N\setminus B_R}\bar{c}e^{-|x|}\le Ce^{-R}R^N\to0,\quad\text{as}\,\,\,R\to+\infty.
\end{equation}
Since $\bar{m}_\e\to\bar{m}$ and $\nabla\bar{u}_\e\to\nabla\bar{u}$ locally uniformly as $\e\to0$, there exists $\e_0$, which depends on $R$, such that $\forall\e\le\e_0$
\begin{equation}\label{5.36}
|\bar{m}_\e-\bar{m}|+\Big||\nabla\bar{u}_\e|^{\g-2}\nabla\bar{u}_\e-\,|\nabla \bar{u}|^{\g-2}\nabla\bar{u}\Big|\le\bar{c}e^{-R},\quad\text{for}\,\,\,|x|\le R.
\end{equation}
Moreover, $\forall\e\le\e_0$ 
$$\bar{m}_\e-\bar{m}+2\nu_R\ge\nu_R, \quad\forall x\in\R^N$$
and hence
\begin{equation}\label{5.38}
|\bar{m}_\e-\bar{m}|\le\bar{m}_\e-\bar{m}+2\nu_R.    
\end{equation}
\indent We estimate each term of the energy $\E_\e$ separately. Concerning the kinetic term, notice that the function $(m,w)\mapsto\frac{m}{\g'}\left|\frac{w}{m} \right|^{\g'}$ is convex and in particular
$$\partial_m\left(\frac{m}{\g'}\left|\frac{w}{m}\right|^{\g'}\right)=-\frac{1}{\g}\left|\frac{w}{m}\right|^{\g'} \qquad\quad\text{and}\qquad\quad \nabla_w\left(\frac{m}{\g'} \left|\frac{w}{m}\right|^{\g'}\right)=\frac{w}{m}\left|\frac{w}{m}\right|^{\g'-2}.$$
By convexity, we estimate separately the integral over $B_R$ and the integral over $\R^N\setminus B_R$, obtaining
\begin{align}\notag
\int_{B_R}\frac{\bar{m}_\e}{\g'}\left|\frac{\bar{w}_\e}{\bar{m}_\e}\right|^{\g'}dx\ge&\int_{B_R}\frac{\bar{m}}{\g'}\left|\frac{\bar{w}}{\bar{m}}\right|^{\g'} dx\\ \label{5.43}
&+\int_{B_R}\frac{\bar{m}_\e-\bar{m}+2\nu_R}{\g'}\left|\frac{\bar{w}_\e-\bar{w}+2\nabla\nu_R}{\bar{m}_\e-\bar{m}+2\nu_R}\right|^{\g'}dx-CR^Ne^{-R},
\end{align}
\begin{align}\notag
\int_{\R^N\setminus B_R}\frac{\bar{m}_\e}{\g'}\left|\frac{\bar{w}_\e} {\bar{m}_\e}\right|^{\g'}dx\ge&\int_{\R^N\setminus B_R}\frac{\bar{m}}{\g'} \left|\frac{\bar{w}}{\bar{m}}\right|^{\g'}dx\\ \label{5.42}
&+\int_{\R^N\setminus B_R}\frac{\bar{m}_\e-\bar{m}+2\nu_R} {\g'}\left|\frac{\bar{w}_\e-\bar{w}+2\nabla\nu_R}{\bar{m}_\e-\bar{m}+2\nu_R}\right|^{\g'}dx-CR^{N+b}e^{-R}.
\end{align}
we also have (refer to estimate (5.44) in \cite{CC}) 
\begin{equation}\label{5.44}
\int_{\R^N}V_\e(x+y_\e)\bar{m}_\e dx\ge\int_{\R^N}V_\e(x+y_\e)\bar{m}dx+ \int_{\R^N}V_\e(x+y_\e)(\bar{m}_\e-\bar{m}+2\nu_R)dx-CR^{b+N}e^{-R}.
\end{equation}
Regarding the Riesz term in the energy $\E_\e$, since by Proposition \ref{prop5.3} $\bar{m}_\e(x)\to\bar{m}(x)$ a.e. as $\e\to0$ and $(\bar{m}_\e)_\e$ is a bounded sequence in $L^\frac{2N}{N+\al}(\R^N)$ (it follows by interpolation using the uniform estimate \eqref{5.16}), applying Lemma \ref{brezis lieb} we get that
\begin{align*} 
\lim\limits_{\e\to0}\int\limits_{\R^N}\int\limits_{\R^N}&\frac{\bar{m}_\e(x)\bar{m}_\e(y)}{|x-y|^{N-\al}}dx\,dy=\int\limits_{\R^N}\int\limits_{\R^N}\frac{\bar{m}(x)\bar{m}(y)}{|x-y|^{N-\al}}dx\,dy+\lim\limits_{\e\to0}\int\limits_{\R^N}\int\limits_{\R^N}\frac{|\bar{m}_\e(x)-\bar{m}(x)|\,|\bar{m}_\e(y)-\bar{m}(y)|}{|x-y|^{N-\al}}dx\,dy\\
&\le\int\limits_{\R^N}\int\limits_{\R^N}\frac{\bar{m}(x)\bar{m}(y)}{|x-y|^{N-\al}}dx\,dy+\lim\limits_{\e\to0}\int\limits_{\R^N}\int\limits_{\R^N}\frac{(\bar{m}_\e-\bar{m}+2\nu_R)(x)\,(\bar{m}_\e-\bar{m}+2\nu_R)(y)}{|x-y|^{N-\al}}dx\,dy
\end{align*}
where in the last inequality we used \eqref{5.38}. Hence
\begin{equation} \label{5.45}
\int_{\R^N}(K_\al\ast\bar{m}_\e)\bar{m}_\e\le
\int_{\R^N}(K_\al\ast\bar{m})\bar{m}+\int_{\R^N}\big(K_\al\ast(\bar{m}_\e-\bar{m}+2\nu_R)\big)(\bar{m}_\e-\bar{m}+2\nu_R)+o_\e(1).
\end{equation}
Finally, putting together estimates \eqref{5.43}, \eqref{5.42}, \eqref{5.44} and \eqref{5.45}, we obtain \eqref{E_e>=}.
\endproof

We are now in position to prove that there is no loss of mass passing to the limit.

\begin{theorem}\label{no loss mass barm}
Let $N-\g'<\al<N$, assume that $(\bar{m}_\e,\bar{w}_\e)$ is a minimizer of $\E_\e$, $\bar{m}$ and $\bar{u}$ obtained from Proposition \ref{prop5.3}. Then,
$$\int_{\R^N}\bar{m}\,dx=M$$
and hence $\bar{m}_\e\to\bar{m}$ in $L^1(\R^N)$. Moreover, for every $\eta>0$, there exist $R,\e_0>0$ such that for all $\e\le\e_0$
\begin{equation}\label{stima mbar B_R}
\int_{B(0,R)}\bar{m}_\e(x)\,dx\ge M-\eta,
\end{equation}
namely 
\begin{equation}\label{stima m Beps}
\int_{B(\e^\frac{\g'}{\g'-N+\al}y_\e,\e^\frac{\g'}{\g'-N+\al}R)}m_\e(x)\,dx\ge M-\eta.
\end{equation}
\end{theorem}

\proof
From \eqref{5.35} we get
$$\int_{\R^N}(\bar{m}_\e-\bar{m}+2\nu_R)dx=M-a+2\int_{\R^N}\nu_Rdx\to M-a, \quad\text{as}\,\,\,R\to+\infty.$$
Let us define
$$C_R:=\frac{M-a}{M-a+2\int_{\R^N}\nu_R}, $$
we observe that $0<C_R<1$ and $C_R\to1$ as $R\to+\infty$, moreover the couple $\big(C_R(\bar{m}_\e- \bar{m}+2\nu_R), C_R(\bar{w}_\e-\bar{w}+2\nabla\nu_R)\big)\in\K_{M-a}$ and it follows that
\begin{align*}
C_R\E_\e&(\bar{m}_\e-\bar{m}+2\nu_R, \bar{w}_\e-\bar{w}+2\nabla\nu_R)= \E_\e\Big(C_R(\bar{m}_\e-\bar{m}+2\nu_R), C_R(\bar{w}_\e-\bar{w}+2\nabla \nu_R)\Big)\\
&+\frac{C_R^2-C_R}{2}\int_{\R^N}\int_{\R^N}\frac{(\bar{m}_\e-\bar{m}+2\nu_R)(x)\,(\bar{m}_\e-\bar{m}+2\nu_R)(y)}{|x-y|^{N-\al}}dx\,dy.
\end{align*}
Notice that $C_R^2-C_R<0$ and from \eqref{5.16} and the fact that $\|\bar{m}_\e\|_{L^1 (\R^N)}=M$, by interpolation we get that $\|\bar{m}_\e\|_{L^p(\R^N)}\le C$ uniformly in $\e$ for every $p\in(1,+\infty)$. By \eqref{disug hls} we get that there exists a constant $C$ independent of $\e$ such that
\begin{align*}
\bigg|\int_{\R^N}\int_{\R^N}&\frac{(\bar{m}_\e-\bar{m}+2\nu_R)(x)\,(\bar{m}_\e-\bar{m}+2\nu_R)(y)}{|x-y|^{N-\al}}dx\,dy\bigg|\le C\|\bar{m}_\e-\bar{m}+2\nu_R\|^2_{L^\frac{2N} {N+\al}(\R^N)}\\
&\le C\left(\|\bar{m}_\e\|_{L^\frac{2N}{N+\al}(\R^N)}+\|\bar{m}\|_{L^\frac{2N}{N+\al}(\R^N)}+2\|\nu_R\|_{L^\frac{2N}{N+\al}(\R^N)}\right)^2\le C.
\end{align*}
Hence
\begin{align*}
C_R\E_\e(\bar{m}_\e-\bar{m}+2\nu_R,\bar{w}_\e-\bar{w}+2\nabla\nu_R)& \ge\E_\e\Big(C_R(\bar{m}_\e-\bar{m}+2\nu_R), C_R(\bar{w}_\e-\bar{w}+2\nabla \nu_R)\Big)
+C\frac{C_R^2-C_R}{2}\\
&\ge \tilde{e}_\e(M-a)+C\frac{C_R^2-C_R}{2}.
\end{align*}
Using this in \eqref{E_e>=} we have that
\begin{align*}
\tilde{e}_\e(M)\ge\tilde{e}_\e(a)+\tilde{e}_\e(M-a)&+(1-C_R)\E_\e(\bar{m}_\e-\bar{m}+2\nu_R,\bar{w}_\e-\bar{w}+2\nabla\nu_R)\\
&+o_\e(1)-CR^{b+N}e^{-R}+C(C_R^2-C_R)
\end{align*}
by \eqref{5.29} we have $\E_\e(\bar{m}_\e-\bar{m}+2\nu_R,\bar{w}_\e- \bar{w}+2\nabla\nu_R)\ge-K$, hence
$$\tilde{e}_\e(M)\ge\,\tilde{e}_\e(a)+\tilde{e}_\e(M-a)-(1-C_R)K+o_\e(1)-CR^{b+N}e^{-R}+C(C_R^2-C_R)$$
finally, from Lemma \ref{lemma5.5} we get
$$0>-C>\,-(1-C_R)K+o_\e(1)-CR^{b+N}e^{-R}+C(C_R^2-C_R)$$
letting $R\to+\infty$ this yields a contradiction. We can conclude following the proof of \cite[Corollary 5.7]{CC}.
\endproof


\subsection{Proof of Theorem \ref{teo no pot}}

We are ready to prove that the triple $(\bar{u}_\e,\bar{m}_\e,\tilde{\la}_\e)$ converges to $(\bar{u},\bar{m},\bar{\la})$ solution to the MFG system \eqref{1 senza V}.

\proof[Proof of Theorem \ref{teo no pot}]

Let $(\bar{u},\bar{m},\bar{\la})$ be the triple obtained from Proposition \ref{prop5.3} and $\bar{w}:=-\bar{m} \nabla\bar{u}|\nabla\bar{u}|^ {\g-2}$.
We have that $(\bar{m},\bar{w})\in\mathcal{B}$, indeed from Proposition \ref{prop5.3} and Theorem \ref{no loss mass barm} we get that $(\bar{m},\bar{w})\in \K_{1,M}$ , and using estimate \eqref{barm<ce^-cx} it follows
$$\int_{\R^N}\bar{m}(1+|x|^b)dx\le\int_{\R^N}c_1e^{-c_2|x|}(1+|x|)^bdx< +\infty.$$
Since $\bar{m}_\e\to\bar{m}$ in $L^1(\R^N)$, $(\bar{m}_\e,\bar{w}_\e),\,(\bar{m}, \bar{w})\in K_{1,M}$ and $\bar{E}_\e,\bar{E}<+\infty$, using Corollary \ref{K_al in Lp e compattezza} we get that $(K_\al\ast\bar{m}_\e)\bar{m}_\e\to(K_\al \ast\bar{m})\bar{m}$ in $L^1(\R^N)$. Moreover, $\bar{w}_\e\to\bar{w}$ locally uniformly and weakly in $L^\frac{\g'\beta}{\g'-1+\beta}(\R^N)$.\\
It follows that the energy $\E_0$ is lower semi-continuous and we have
\begin{align}\notag
\E_0(\bar{m},\bar{w})&=\int_{\R^N}\frac{\bar{m}}{\g'}\left|\frac{\bar{w}}{\bar{m}}\right|^{\g'}dy-\int_{\R^N}\int_{\R^N}\frac{\bar{m}(x)\,\bar{m}(y)}{|x-y|^{N-\al}}dx\,dy\\\label{gamma}
&\le\liminf\limits_{\e}\left(\int_{\R^N}\frac{\bar{m}_\e}{\g'}\left|\frac{\bar{w}_\e}{\bar{m}_\e}\right|^{\g'}dy-\int_{\R^N}\int_{\R^N}\frac{\bar{m}_\e(x)\,\bar{m}_\e(y)}{|x-y|^{N-\al}}dx\,dy\right)\le\liminf\limits_{\e}\E_\e(\bar{m}_\e,\bar{w}_\e) 
\end{align}
where in the last inequality we used the fact that $V\ge0$. 
Moreover, if $(m,w)\in\mathcal{B}$, using \eqref{V_e}, we have
$$0\le\lim\limits_{\e\to0}\int_{\R^N}m(y+y_\e)V_\e(y+y_\e)dy\le\lim\limits_{\e\to0}C_V\e^{\frac{(N-\al)\g'}{\g'-N+\al}}\int _{\R^N}(1+|y|)^b m(y)dy=0$$
from which it follows that
$$\lim\limits_{\e\to0}\E_\e\big(m(\cdot+y_\e), w(\cdot+y_\e)\,\big)=\E_0(m,w).$$
Using the fact that $(\bar{m}_\e,\bar{w}_\e)$ is a minimizer for $\E_\e$ and then \eqref{gamma}, we finally get
$$\E_0(m,w)=\lim\limits_{\e\to0}\E_\e\big(m(\cdot+y_\e), w(\cdot+y_\e)\,\big)\ge\lim\limits_{\e\to0}\E_\e(\bar{m}_\e,\bar{w}_\e)\ge\E_0(\bar{m},\bar{w}),$$
this proves that 
$$\E_0(\bar{m},\bar{w})=\min\limits_{(m,w)\in\mathcal{B}}\E_0(m,w).$$
Since $(\bar{m}_\e,\bar{w}_\e)$ and $(\bar{m},\bar{w})$ are minimizers of $\E_\e$ and $\E_0$ respectively, we obtain that
\begin{align*}
\E_\e(\bar{m}_\e,\bar{w}_\e)&\le\int_{\R^N}\frac{\bar{m}(y+y_\e)}{\g'}\left|\frac{\bar{w}(y+y_\e)}{\bar{m}(y+y_\e)}\right|^{\g'}+V_\e(y+y_\e)\,\bar{m}(y+y_\e)-\bar{m}(y+y_\e) K_\al\ast\bar{m}(y+y_\e)dy\\
&=\E_0(\bar{m},\bar{w})+\int_{\R^N}V_\e(y+y_\e)\,\bar{m}(y+y_\e)dy\le 
\E_0(\bar{m}_\e,\bar{w}_\e)+C\e^{\frac{(N-\al)\g'}{\g'-N+\al}}
\end{align*}
where in the last inequality we used also \eqref{V_e} and the fact that $(\bar{m},\bar{w})\in\mathcal{B}$. It follow immediately that 
$$\int_{B(0,R)}\bar{m}_\e(y)\,V_\e(y+y_\e)dy\le C\e^{\frac{(N-\al)\g'} {\g'-N+\al}}.$$
From \eqref{V_e} and \eqref{stima mbar B_R} we get that there exists a positive constant $C$ such that
\begin{equation}\label{eps y_e<C}
\e^{\frac{\g'}{\g'-N+\al}}|y_\e|\le C
\end{equation}
and hence using \eqref{V_e} again we obtain that
$$0\le V_\e(y+y_\e)\le C_V\e^{\frac{(N-\al)\g'}{\g'-N+\al}}\left(1+\e^\frac{\g'} {\g'-N+\al}|y+y_\e|\right)^b\le C\e^{\frac{(N-\al)\g'}{\g'-N+\al}} \left(1+\e^\frac{\g'}{\g'-N+\al}|y|\right)^b$$
which allows us to conclude that $V_\e(y+y_\e)\to0$ locally uniformly as $\e\to0$. This prove that the function $g$ defined in Theorem \ref{prop5.3} is actually zero.
\endproof


\section{Concentration of the mass}

The following result allows us to localize the points where the mass concentrates. 

\begin{prop}\label{prop localize point conc}
As $\e\to0$, the sequence $\e^\frac{\g'}{\g'-N+\al}y_\e$ converges, up to subsequences, to a point $\bar{x}\in\R^N$ such that $V(\bar{x})=0$.
\end{prop}

\proof
Let $z\in\R^N$ (to be fixed later), by minimality of $(\bar{m}_\e,\bar{w}_\e)$ and $(\bar{m}(\cdot+z),\bar{w}(\cdot+z))$ we get
\begin{align*}
    \E_\e(\bar{m}_\e,\bar{w}_\e)&\le\E_\e\big(\bar{m}(\cdot+z),\bar{w}(\cdot+z)\big)=\E_0(\bar{m},\bar{w})+\int_{\R^N}\bar{m}(y+z)\,V_\e(y+y_\e)dy\\
    &\le\int_{\R^N}\frac{\bar{m}_\e}{\g'}\left|\frac{\bar{w}_\e}{\bar{m}_\e}\right|^{\g'}dy-\int_{\R^N}\int_{\R^N}\frac{\bar{m}_\e(x)\,\bar{m}_\e(y)}{|x-y|^{N-\al}}dx\,dy+\int_{\R^N}\bar{m}(y+z)\,V_\e(y+y_\e)dy
\end{align*}
hence
$$\int_{\R^N}\bar{m}_\e(y)\,V_\e(y+y_\e)dy\le\int_{\R^N}\bar{m}(y+z)\,V_\e(y+y_\e)dy=\int_{\R^N}\bar{m}(y)\,V_\e(y+y_\e-z)dy$$
and using \eqref{def V_e} and the fact that $\bar{m}_\e(y)=\e^\frac{N\g'} {\g'-N+\al}m_\e\big(\e^\frac{\g'}{\g'-N+\al}y+\e^\frac{\g'}{\g'-N+\al}y_\e\big)$
we get
$$\int_{\R^N}m_\e(y)\,V(y)dy\le\int_{\R^N}\bar{m}(y)\,V\big(\e^\frac{\g'}{\g'-N+\al}(y+y_\e-z)\big)dy.$$
By assumption the potential $V$ is a locally H\"older continuous coercive function, so it has a global minimum at a point $\bar{z}\in\R^N$, and by a shift of $\la$ we may assume that $V(\bar{z})=0$. Let us fix $z=y_\e-\e^{-\frac{\g'}{\g'-N+\al}}\bar{z}$, it holds
\begin{equation}\label{limsup int mV}
\limsup\limits_{\e\to0}\int_{\R^N}\bar{m}(y)\,V\big(\e^\frac{\g'}{\g'-N+\al}y+\bar{z}\big)dy\le \limsup\limits_{\e\to0}c_1\int_{\R^N}e^{-c_2|y|}V\big(\e^\frac{\g'}{\g'-N+\al}y+\bar{z}\big)dy=0.  
\end{equation}
Moreover, by \eqref{eps y_e<C}, we get that (up to subsequences) 
$$\e^{\frac{\g'}{\g'-N+\al}}y_\e\to\bar{x}\in\R^N$$
and by \eqref{stima m Beps} denoting by $B:=B\big(\e^\frac{\g'}{\g'-N+\al}y_\e,\e^\frac{\g'}{\g'-N+\al}R\big)$
\begin{align}\label{liminf int mV}
\liminf\limits_{\e\to0}&\int_{\R^N}m_\e(y)V(y)dy\ge\liminf\limits_{\e\to0}\int_{B}m_\e(y)V(y)dy\\
&\ge\liminf\limits_{\e\to0}\,\min\limits_{x\in B}V(x)\int_{B}m_\e(y)dy\ge\liminf \limits_{\e\to0}\,\min\limits_{x\in B}V(x)(M-\eta)\\
&\ge(M-\eta)V(\bar{x}).
\end{align}
From \eqref{limsup int mV} and \eqref{liminf int mV} we can obtain that $V(\bar{x})=0$.
\endproof

\proof[Proof of Theorem \ref{TEO concentrazione massa}] It follows from Proposition \ref{prop localize point conc} and Theorem \ref{no loss mass barm}.
\endproof

\begin{obs}
Arguing as in \cite{CC} (refer to Proposition 5.13 and its proof) one can prove that if $V$ has a finite number of minima $x_i\in\R^N$ for $i=1,\dots,n$ and can be written as 
$$V(x)=h(x)\prod\limits_{i=1}^n|x-x_i|^{b_i}$$
for a certain function $C_V^{-1}\le h(x)\le C_V$ on $\R^N$ and $b_i>0$ such that $\sum_{i=1}^n b_i=b$, then the sequence $\e^\frac{\g'}{\g'-N+\al}y_\e$, as $\e\to0$, converges (up to subsequences) to the more stable minimum of $V$ (namely the point $x_j$ such that $b_j=\max\limits_{i=1...n}b_i$).
\end{obs}

\textbf{Acknowledgments.} The author would like to thank Professor Annalisa Cesaroni for helpful suggestions and fruitful discussions.


\end{document}